\documentclass[11pt]{amsart}

\usepackage{amsaddr}

\usepackage[misc,geometry]{ifsym}

\usepackage[english]{babel} 
\usepackage[UKenglish]{datetime}

\usepackage{ae, aecompl}
\usepackage{amscd}
\usepackage{amsfonts}
\usepackage{amsmath}
\usepackage{amssymb}
\usepackage{amsthm}
\usepackage{amsxtra}
\usepackage{cases}
\usepackage{chngcntr}
\usepackage{cite}
\usepackage{color}
\usepackage{graphicx}
\usepackage{latexsym}
\usepackage{bbm} 
\usepackage{bm}
\usepackage{bbold}
\usepackage{mathtools}
\usepackage{microtype}
\usepackage{qsymbols}
\usepackage[active]{srcltx}
\usepackage{tikz}
\usepackage{tikz-cd}
\usepackage{url}
\usepackage{verbatim} 
\usepackage[all]{xy}

\usepackage{hhline}
\usepackage{array}
\usepackage{longtable}

\usepackage{hyperref}
\usepackage{enumitem}
\usepackage{cleveref}

\tikzcdset{arrow style=tikz,
	diagrams={>={Straight Barb[scale=0.7]}}
}

\usepackage{tikz,tikz-cd}
\usetikzlibrary{matrix, arrows}

\makeatletter
\tikzset{
	edge node/.code={%
		\expandafter\def\expandafter\tikz@tonodes\expandafter{\tikz@tonodes #1}}}
\makeatother
\tikzset{
	subset/.style={
		draw=none,
		edge node={node [sloped, allow upside down, auto=false]{$\subset$}}},
	subseteq/.style={
		draw=none,
		every to/.append style={
			edge node={node [sloped, allow upside down, auto=false]{$\subseteq$}}}}
}

\usepackage[charter]{mathdesign}
\newcommand{\bbfont}{\mathbbm}
\newcommand{\catfont}[1]{\textbf{\fontfamily{qag}\selectfont{#1}}}
\newcommand{\catfontscript}[1]{\textbf{\fontfamily{qag}\selectfont{\footnotesize{#1}}}}

\makeatletter
\newcommand{\btfs}[1]
{
	\ifthenelse{\equal{\f@shape}{n}}{\ensuremath{\mathrm{#1}}}
	{\ifthenelse{\equal{\f@shape}{sc}}{\ensuremath{\mathrm{#1}}}
		{\ifthenelse{\equal{\f@shape}{it}}{\ensuremath{\mathit{#1}}}
			{\ifthenelse{\equal{\f@shape}{sl}}{\ensuremath{\mathit{#1}}}{}	
			}
		}
	}
}
\makeatother

\renewcommand{\btfs}{\mathrm}

\newcommand{\NN}{{\bbfont N}}

\newcommand{\abs}[1]{{\lvert #1 \rvert}}
\newcommand{\norm}[1]{{\lVert #1 \rVert}}
\newcommand{\lrnorm}[1]{{\left\lVert #1 \right\rVert}}
\newcommand{\lrabs}[1]{{\left\lvert #1 \right\rvert}}

\newcommand{\cont}{\text{C}}
\newcommand{\contc}{\cont_{\text{c}}}

\theoremstyle{plain}

\newtheorem{theorem}{Theorem}[section]
\newtheorem{proposition}[theorem]{Proposition}
\newtheorem{lemma}[theorem]{Lemma}
\newtheorem{corollary}[theorem]{Corollary}

\newtheorem*{theorem*}{Theorem}
\newtheorem*{proposition*}{Proposition}
\newtheorem*{lemma*}{Lemma}
\newtheorem*{corollary*}{Corollary}

\theoremstyle{definition}

\newtheorem{example}[theorem]{Example}
\newtheorem{remark}[theorem]{Remark}

\newtheorem*{definition*}{Definition}
\newtheorem*{example*}{Example}
\newtheorem*{remark*}{Remark}

\setlist[enumerate,1]{label=\textup{(\arabic*)},ref=(\arabic*)}
\setlist[enumerate,2]{label=\textup{(\alph*)},ref=(\arabic{enumi}.\alph*)}
\setlist[enumerate,3]{label=\textup{(\roman*)},ref=(\arabic{enumi}.\alph{enumii}.\roman*)}
\setlist[enumerate,4]{label=\textup{(\Alph*)},ref=(\arabic{enumi}.\alph{enumii}.\roman{enumiii}.\Alph*)}

\crefname{theorem}{Theorem}{Theorems}
\crefname{proposition}{Proposition}{Propositions}
\crefname{lemma}{Lemma}{Lemmas}
\crefname{corollary}{Corollary}{Corollaries}
\crefname{definition}{Definition}{Definitions}
\crefname{example}{Example}{Examples}
\crefname{remark}{Remark}{Remarks}

\crefname{section}{Section}{Sections}
\crefname{subsection}{Section}{Sections}
\crefname{subsubsection}{Section}{Sections}

\crefname{equation}{equation}{equations}

\numberwithin{equation}{section}

\allowdisplaybreaks 

\newcommand{\dc}{\text{d}}

\newcommand{\map}{\varphi}
\newcommand{\maptwo}{\psi}
\newcommand{\ndual}[1]{#1^\ast}
\newcommand{\nbidual}[1]{#1^{\ast\ast}}
\newcommand{\odual}[1]{#1^\sim}

\newcommand{\pos}[1]{#1^+}

\newcommand{\veps}{\varepsilon}
\newcommand{\supp}{\operatorname{supp}}

\newcommand{\Ell}{\text{L}}

\newcommand{\cm}{\varphi} 
\newcommand{\mil}{\varphi} 
\newcommand{\indexset}{I}
\newcommand{\indexone}{i}
\newcommand{\indextwo}{j}
\newcommand{\obj}{O}
\newcommand{\factor}{\chi}

\newcommand{\systemofobjects}[3]{(#1_{#2})_{#3}}
\newcommand{\systemofconnectingmorphisms}[3]{(#1_{#2})_{#3}}
\newcommand{\systemofmorphisms}[3]{(#1_{#2})_{#3}}
\newcommand{\directedsystem}[6]{\bigl(\systemofobjects{#1}{#2}{#3},\systemofconnectingmorphisms{#4}{#5}{#6}\big)}

\newcommand{\Odirsysfull}{\directedsystem{\obj}{\indexone}{\indexone\in\indexset}{\cm}{\indextwo\indexone}{\indexone,\indextwo\in\indexset,\indextwo\geq\indexone}}
\newcommand{\Odirsys}{\directedsystem{\obj}{\indexone}{}{\cm}{\indextwo\indexone}{\indextwo\geq\indexone}}
\newcommand{\Edirsys}{\directedsystem{E}{\indexone}{}{\cm}{\indextwo\indexone}{\indextwo\geq\indexone}}

\newcommand{\Osys}{\bigl(\obj,\systemofmorphisms{\mil}{\indexone}{}\big)}
\newcommand{\Osysprimed}{\bigl(\obj^\prime,\systemofmorphisms{\mil^\prime}{\indexone}{}\big)}

\newcommand{\Esys}{\bigl(E,\systemofmorphisms{\mil}{\indexone}{}\big)}

\newcommand{\Esysprimed}{\bigl(E^\prime,\systemofmorphisms{\mil^\prime}{\indexone}{}\big)}

\newcommand{\CAT}{\catfont{Cat}}

\newcommand{\SET}{\catfont{Set}}
\newcommand{\MET}{\catfont{Met}}
\newcommand{\COMMET}{\catfont{ComMet}}

\newcommand{\POVSP}{\catfont{POVS}_{\catfontscript{Pos}}}
\newcommand{\VS}{\catfont{VS}_{\catfontscript{}}}

\newcommand{\PONSP}{\catfont{PONS}_{\catfontscript{Pos}}}
\newcommand{\NS}{\catfont{NS}_{\catfontscript{}}}

\newcommand{\POBSP}{\catfont{POBS}_{\catfontscript{Pos}}}
\newcommand{\BS}{\catfont{BS}_{\catfontscript{}}}

\newcommand{\VLIPLH}{\catfont{VL}_{\catfontscript{IPLH}}}
\newcommand{\VLLH}{\catfont{VL}_{\catfontscript{LH}}}
\newcommand{\VLIP}{\catfont{VL}_{\catfontscript{IP}}} 

\newcommand{\NLAIPLH}{\catfont{NL}_{\catfontscript{AIPLH}}}
\newcommand{\NLIPLH}{\catfont{NL}_{\catfontscript{IPLH}}}
\newcommand{\NLLH}{\catfont{NL}_{\catfontscript{LH}}}
\newcommand{\NLAIP}{\catfont{NL}_{\catfontscript{AIP}}} 
\newcommand{\NLIP}{\catfont{NL}_{\catfontscript{IP}}} 

\newcommand{\BLAIPLH}{\catfont{BL}_{\catfontscript{AIPLH}}}
\newcommand{\BLIPLH}{\catfont{BL}_{\catfontscript{IPLH}}}
\newcommand{\BLLH}{\catfont{BL}_{\catfontscript{LH}}}
\newcommand{\BLAIP}{\catfont{BL}_{\catfontscript{AIP}}} 
\newcommand{\BLIP}{\catfont{BL}_{\catfontscript{IP}}} 

\begin{document}


\title[Direct limits]{Direct limits in categories of normed vector lattices and Banach lattices}

\author{\Letter\  Chun Ding}

\address[\Letter Chun Ding]{Mathematical Institute, Leiden University, P.O.\ Box 9512, 2300 RA Leiden, the Netherlands\\ $\text{c}\_{}\text{ding@foxmail.com}$}

\author{Marcel de Jeu}
\address[Marcel de Jeu]{Mathematical Institute, Leiden University, P.O.\ Box 9512, 2300 RA Leiden, the Netherlands\\ and\\
	Department of Mathematics and Applied Mathematics, University of Pretoria, Corner of Lynnwood Road and Roper Street, Hatfield 0083, Pretoria,
	South Africa\\
mdejeu@math.leidenuniv.nl}




\subjclass[2010]{Primary 46M40; Secondary 46A40, 46B40, 46B42, 46E05}

\keywords{Vector lattice, normed vector lattice, Banach lattice, direct limit, inductive limit, (almost) interval preserving map, order continuous norm}

\begin{abstract}
After collecting a number of results on interval and almost interval preserving linear maps and vector lattice homomorphisms, we show that direct systems in various categories of normed vector lattices and Banach lattices have direct limits, and that these coincide with direct limits of the systems in naturally associated other categories.  For those categories where the general constructions do not work to establish the existence of general direct limits, we describe the basic structure of those direct limits that \emph{do} exist.\\
A direct system in the category of Banach lattices and contractive almost interval preserving vector lattice homomorphisms has a direct limit. When the Banach lattices in the system all have order continuous norms, then so does the Banach lattice in a direct limit. This is used to show that a Banach function space over a locally compact Hausdorff space has an order continuous norm when the topologies on all compact subsets are metrisable and (the images of) the continuous compactly supported functions are dense.
\end{abstract}

\maketitle



\section{Introduction and overview}\label{sec:introduction_and_overview}


\noindent Until very recently, there has been a rather modest role for direct limits of vector lattices in the literature.
In \cite{shannon:1977}, Shannon characterises the spaces of continuous compactly supported functions on locally compact Hausdorff spaces as the direct limits of certain direct systems of Banach lattices in the category of normed vector lattices and isometric lattice homomorphisms. In \cite{filter:1988}, Filter shows that direct limits exist in the category of vector lattices and lattice homomorphisms.  He also studies permanence properties when the lattice homomorphisms are injective.
It is only recently in \cite{van_amstel_van_der_walt_UNPUBLISHED:2022} that the subject has been taken up again, and more comprehensively, by Van Amstel and Van der Walt. They extend the work in \cite{filter:1988} on direct limits of vector lattices by also taking interval preserving and order continuous lattice homomorphisms into account. In addition, they set up the basic theory for inverse limits of vector lattices, for which there is also only a very limited literature (see \cite{van_amstel_van_der_walt_UNPUBLISHED:2022} for references). Both types of limits are then used to study the relation between a vector lattice and its order dual, and a number of applications of the theory of direct and inverse limits to several problems in concrete vector lattices are given.

The current paper is primarily concerned with the existence of direct limits of direct systems in various categories of normed vector lattices and Banach lattices. The presence of a norm makes this more complicated than the purely algebraic case of vector lattices, where we have only little to add to \cite{van_amstel_van_der_walt_UNPUBLISHED:2022}. Nevertheless, we still also give a method to construct direct limits of direct systems in categories of vector lattices that is different from that in \cite{filter:1988} (as recapitulated in \cite{van_amstel_van_der_walt_UNPUBLISHED:2022}). It appears to be a somewhat more transparent and, what is more, it is naturally modified to show that direct limits also exist in various categories of normed vector lattices and Banach lattices and contractive vector lattice homomorphisms.

\smallskip

\noindent This paper is organised as follows.

\cref{sec:preliminaries} contains preliminary material. The definitions of a direct limit in a category and of (almost) interval preserving linear maps are recalled, and the thirteen categories of vector lattices, normed vector lattices, and Banach lattices are introduced that we shall consider in this paper, together with a few other ones that occur naturally.

\cref{sec:interval_preserving_and_almost_interval_preserving_linear_maps_and_lattice_homomorphisms} on interval and almost interval preserving linear maps and lattice homomorphisms is the toolbox for this paper. We have collected the basic general results in the literature about these maps that we are aware of, and added our own contributions. We believe that this section contains a fairly comprehensive basic theory of the categorical aspects of such maps that could perhaps also find use elsewhere.

In \cref{sec:direct_limits_existence_via_three_standard_constructions}, we give what we call the standard constructions of direct limits of (suitable) direct systems of\textemdash in order of increasing complexity of the constructions\textemdash vector lattices, normed vector lattices, and Banach lattices. They are essentially the well-known basic constructions of direct limits of direct systems of vector spaces and linear maps, and of normed spaces and Banach spaces and contractions, respectively. To see why these also work for some categories of vector lattices, normed vector lattices, and Banach lattices, and not for others, it is necessary to briefly go through the three basic constructions. These reviews of the constructions also make it immediately clear why the direct limits of a fixed direct system in various categories coincide. We also show that a direct limit of a direct system of normed spaces (resp.\ Banach spaces) and contractive linear maps is also a direct limit of that system in the category of metric spaces (resp.\ complete metric spaces) and contractive maps. Although we may not be the first to note this, we are not aware of a reference for these two facts.

The constructions in \cref{sec:direct_limits_existence_via_three_standard_constructions} are not guaranteed to work for all thirteen categories of vector lattices, normed vector lattices, and Banach lattices under consideration. For six of these, the existence of direct limits of their general direct systems remains open. Nevertheless, the results in \cref{sec:interval_preserving_and_almost_interval_preserving_linear_maps_and_lattice_homomorphisms}  enable us to say something in \cref{sec:direct_limits_additional_results} about the basic structure of those direct limits that \emph{do} exist. For five of these exceptional categories, we give an example of a direct system for which the standard construction `unexpectedly' still produces a direct limit, and an example where it fails to do so.

The concluding \cref{sec:direct_limits_and_order_continuity} is concerned with the order continuity of direct limits of direct systems of Banach lattices. \cref{sec:direct_limits_existence_via_three_standard_constructions} shows that direct limits exist in the category of Banach lattices and contractive almost interval preserving vector lattice homomorphisms. When all Banach lattices in a direct system in this category are order continuous, then, as is shown in this section, so is the Banach lattice in its direct limit. As an application, we show that a Banach function space over a locally compact Hausdorff space has an order continuous norm when the topologies on all its compact subsets are metrisable and (the images of) the continuous compactly supported functions are dense in it. An alternate, more direct proof of this is also provided.

\smallskip

\noindent It seems natural to wonder whether methods as in \cite{van_amstel_van_der_walt_UNPUBLISHED:2022}, using direct as well as inverse limits in categories of vector lattices, can also be applied to problems in normed vector lattices and Banach lattices, once the material in the present paper has been supplemented with sufficiently many results on inverse limits of direct systems of normed vector lattices and Banach lattices.


\section{Preliminaries}\label{sec:preliminaries}


\noindent In this section, we give the conventions, notations, and definitions used in the sequel.

All vector spaces in this paper are real vector spaces.\footnote{The constructions of direct limits of direct systems of vector spaces, normed spaces, and Banach spaces that we shall review also work for complex spaces.}. A preordered vector space is a vector space with a linear preorder induced by a wedge that need not be a cone.  The positive wedges of preordered vector spaces need not be generating.  Vector lattices need not be Archimedean. If $S$ is a subset of a vector lattice~$E$, then we write~$\pos{S}$ for the set $\{\pos{s}: s\in S\}$ of the positive parts of its elements. Hence $\pos{S}\supseteq S\cap \pos{E}$ with a possibly proper inclusion; for a linear subspace~$L$ of~$E$, we have $\pos{L}=L\cap \pos{E}$ if and only if~$L$ is a vector sublattice of~$E$. We write~$\odual{E}$ for the order dual of a vector lattice~$E$. By a lattice homomorphism we mean a vector lattice homomorphism; by a normed lattice we mean a normed vector lattice; and by an order continuous Banach lattice we mean a Banach lattice that has an order continuous vector norm.  The positive wedges of preordered normed spaces need not be closed. A contraction between two normed spaces is supposed to be linear. We write~$\ndual{E}$ for the norm dual of a normed space~$E$.

When~$E$ is a vector lattice and $e\leq e^\prime$ in~$E$, then we let $[e,e^\prime]_E\coloneqq\{e^{\prime\prime}\in E: e\leq e^{\prime\prime} \leq e^\prime\}$ denote the corresponding order interval in~$E$.

A linear map $\map\colon E\to F$ between two vector lattices is said to be \emph{interval preserving} if it is positive and such that $\map([0,x]_E)=[0,\map(x)]_F$ for all $x\in\pos{E}$.

A linear map $\map\colon E\to F$ from a vector lattice~$E$ into a normed lattice~$F$ is called \emph{almost interval preserving} if it is positive and such that $[0, \map(x)]_F=\overline{\map([0,x]_E)}$ for all $x\in\pos{E}$. It is clear that $[0, \map(x)]_F\supseteq\overline{\map([0,x]_E)}$ for positive~$\map$; the point is the reverse inclusion.

Let  $\CAT$ be a category. Suppose that~$\indexset$ is a directed non-empty set, and that $\Odirsysfull$ is a pair, consisting of a collection of objects~$\obj_i$ indexed by~$\indexset$, and morphisms $\cm_{ji}\colon\obj_i\to \obj_j$ for all $i,j\in \indexset$ with  $i\geq j$, such that~$\cm_{ii}$ is the identity morphism of~$\obj_i$ and $\cm_{kj}\circ\cm_{ji}=\cm_{ki}$ for all $i,j,k\in I$ with $i\leq j\leq k$. Then  $\Odirsysfull$ is called a \emph{direct system in $\CAT$} over~$\indexset$. As we shall always be working with one direct system at a time, we shall omit the mention of the index set~$\indexset$ in the proofs and in the notation altogether, and simply write $\Odirsys$. If~$\obj$ is an object in $\CAT$ and $\mil_i\colon\obj_i\to\obj$ are morphisms, then the system $\Osys$ is called \emph{compatible with $\Odirsys$} when $\mil_j\circ\cm_{ji}=\mil_i$ for all $j\geq i$.
A \emph{direct limit} of the system $\Odirsys$ in $\CAT$ is a compatible system $\Osys$ with the property that, for every compatible system  $\Osysprimed$,  there is a unique morphism $\factor\colon\obj\to\obj^\prime$ such that $\mil_i^\prime=\factor\circ\mil_i$ for all~$i$.\footnote{The terminology in the literature is not uniform: `inductive system' and `inductive limit' are also used, for example. In categorical language, a direct limit as above is a co-limit of the diagram that is provided by $\Odirsys$.} The commutativity of the diagram
\[
	\begin{tikzcd}
		\obj_i\arrow{rr}{\cm_{ji}}\arrow{rd}[swap]{\mil_i}\arrow[bend right]{ddr}[swap]{\mil_i^\prime} & & \obj_j\arrow{ld}{\mil_j}\arrow[bend left]{ddl}{\mil_j^\prime}\\
		& \obj\arrow{d}{\factor} &\\
		& \obj^\prime&
	\end{tikzcd}
\]
shows that, in retrospect, the compatibility of $\Osysprimed$ `originates' from the compatibility of $\Osys$.

Direct limits need not exist, but if they do then they are unique up to isomorphism in a strong sense. If $\Osys$ and $\Osysprimed$ are both direct limits of $\Odirsys$, and $\factor\colon\obj\to\obj^\prime$ is the unique morphism such that $\mil_i^\prime=\factor\circ\mil_i$ for all~$i$, then~$\factor$ is an isomorphism. Its inverse is the unique morphism $\factor^\prime\colon\obj^\prime\to\obj$ such that $\mil_i=\factor^\prime\circ\mil_i^\prime$ for all~$i$.

For a category of vector lattices and linear maps, there are two properties that we consider for the linear maps that are its morphisms: being interval preserving and being a lattice homomorphism. This gives four categories of vector lattices, and we ask for the existence of direct limits of direct systems in each of these. One of these categories consists of the vector lattices and the linear maps. In this category, all direct systems have direct limits that coincide with their direct limits in the category of vector spaces and linear maps. Indeed, if~$E$ is the vector space in a direct limit in the latter category, then one merely needs to supply it with the structure of a vector lattice. This is clearly possible since every vector space is isomorphic, as a vector space, to the vector lattice of finitely supported functions on a basis. This leaves three categories to consider.

For a category of normed lattices or Banach lattices and contractions, there are three properties that we consider for the contractions:  being interval preserving, being almost interval preserving, and being a lattice homomorphism. This leads to eight categories for each. Here, we have nothing to say about the (non-trivial) cases where the morphisms are simply contractions, and we shall ignore this case in the sequel. This leaves seven categories for each, but since being interval preserving implies being almost interval preserving there are effectively five categories of normed lattices or Banach lattices to consider in the sequel.

We give an overview of these thirteen categories in the following table, and also include additional categories that will turn out to occur naturally in the sequel.

\bigskip
\medskip

\begin{center}
	\begin{longtable}{ |m{1.7cm}|m{4.9cm}|m{4.9cm}| }
		\hline
		& \sc{Objects}&\sc{Morphisms}\\
		\hhline{|=|=|=|}
		$\VLIPLH$ & Vector lattices &Interval preserving lattice homomorphisms \\
		\hline
		$\VLLH$ & Vector lattices & Lattice homomorphisms \\
		\hline
		$\POVSP$ & Preordered vector spaces & Positive linear maps \\
		\hline
		$\VS$& Vector spaces & Linear maps\\
		\hline
		$\SET$& Sets & Arbitrary maps\\
		\hline
		$\VLIP$ & Vector lattices & Interval preserving linear maps\\
		\hhline{|=|=|=|}
		$\NLIPLH$ & Normed lattices & Contractive interval preserving lattice homomorphisms\\
		\hline
		$\NLAIPLH$ & Normed lattices &Contractive almost interval preserving lattice homomorphisms \\
		\hline
		$\NLLH$ & Normed lattices & Contractive lattice homomorphisms \\
		\hline
		$\PONSP$ & Preordered normed spa\-ces  & Positive contractions\\
		\hline
		$\NS$& Normed spaces & Contractions\\
		\hline
		$\MET$ & Metric spaces & Contractive maps\\
		\hline	
		$\NLIP$ & Normed lattices & Interval preserving contractions\\
		\hline
		$\NLAIP$ & Normed lattices & Almost interval preserving contractions\\
		\hhline{|=|=|=|}
		$\BLAIPLH$ & Banach lattices & Contractive almost interval preserving lattice homomorphisms\\
		\hline
		$\BLLH$ & Banach lattices & Contractive lattice homomorphisms \\
		\hline
		$\POBSP$ & Preordered Banach spaces  & Positive contractions\\
		\hline
		$\BS$ & Banach spaces & Contractions\\
		\hline
		$\COMMET$ & Complete metric spaces & Contractive maps\\
		\hline
		$\BLIP$ & Banach lattices & Interval preserving contractions\\
		\hline
		$\BLAIP$ & Banach lattices & Almost interval preserving contractions\\	
		\hline
		$\BLIPLH$ & Banach lattices &Contractive interval preserving lattice homomorphisms \\
		\hline
	\end{longtable}
\end{center}

The ordering in the table may look a bit odd at first sight, but it reflects the existence of the following chains, in which the inclusion symbols denote subcategory/supercategory relations, and arrows indicate the obvious forgetful functors: 
\begin{align*}
	&&\phantom{\hskip -1em \subset \hskip 1em}&\VLIPLH &\hskip -1em \subset \hskip 1em &\VLLH &\hskip -1em \subset \hskip 1em &\POVSP &\hskip -1em \to \hskip 1em &\VS &\hskip -1em \to \hskip 1em&\SET\\
	&\NLIPLH&\hskip -1em \subset \hskip 1em &\NLAIPLH &\hskip -1em \subset \hskip 1em &\NLLH &\hskip -1em \subset \hskip 1em &\PONSP &\hskip -1em \to \hskip 1em &\NS&\hskip -1em \subset \hskip 1em&\MET\\
	&&\phantom{\hskip -1em \subset \hskip 1em}&\BLAIPLH &\hskip -1em \subset \hskip 1em  &\BLLH &\hskip -1em \subset \hskip 1em &\POBSP &\hskip -1em \to \hskip 1em& \BS&\hskip -1em \subset \hskip 1em&\COMMET
\end{align*}
As we shall see in \cref{res:direct_limits_in_chain_of_vector_lattices}, \cref{res:direct_limits_in_chain_of_normed_vector_lattices}, and \cref{res:direct_limits_in_chain_of_banach_lattices}, respectively, direct systems in those categories in the above chains that consist of vector lattices have direct limits in the pertinent categories, and these direct limits are also direct limits of these systems in the categories to the right in the same chain.\footnote{Strictly speaking, one should say that the images of the direct limits under the appropriate combinations of inclusion and/or forgetful functors are direct limits of the images of the direct systems under the same combinations of functors. Since there appears to be little chance of confusion, we prefer to use a shorter formulation as in the text.} Direct limits also exist in $\VS$, $\NS$, and $\BS$, and these are also direct limits of the system in $\SET$, $\MET$, and $\COMMET$, respectively.

The categories $\VLIP$, $\NLIP$, $\NLAIP$, $\BLIP$, $\BLAIP$, and $\BLIPLH$ are outliers in the sense that the general constructions of direct limits that work for the other categories of vector lattices may fail in these cases. The existence of direct limits of general direct systems in these categories is unclear. We shall comment further on this in \cref{sec:direct_limits_additional_results}.


\section{Interval preserving and almost interval preserving linear maps and lattice homomorphisms}\label{sec:interval_preserving_and_almost_interval_preserving_linear_maps_and_lattice_homomorphisms}


\noindent In this section, we collect a number of results of a categorical flavour on interval and almost interval preserving linear maps and lattice homomorphisms. Our primary goal is an application of such results in the context of direct limits in  \cref{sec:direct_limits_existence_via_three_standard_constructions}. We include more than is needed for just that, however, in an attempt to fill a reasonably complete toolbox that may also find use elsewhere. Some of the results are elementary and have only been included explicitly to complete the picture and for reference purposes in the sequel of the paper, but others are less obvious.

We start with the case of interval preserving maps and vector lattices. After that, we consider that of normed lattices and almost interval preserving maps, where the proofs of analogous and additional results are less straightforward.


\subsection{Interval preserving linear maps and lattice homomorphisms}\label{subsec:interval_preserving_linear_maps_and_lattice_homomorphisms}


Interval preserving linear maps are related to lattice homomorphisms via duality. When $\map\colon E\to F$ is an interval preserving linear map between two vector lattices, then its order adjoint $\odual{\map}\colon\odual{F}\to\odual{E}$ is a lattice homomorphism; when~$\map$ is a lattice homomorphism, then  $\odual{\map}$ is interval preserving; when $\odual{\map}$ is interval preserving and~$\odual{F}$ separates the points of~$F$, then~$\map$ is a lattice homomorphism. We refer to \cite[Theorems~2.16 and 2.20]{aliprantis_burkinshaw_POSITIVE_OPERATORS_SPRINGER_REPRINT:2006} for these results.

The following two lemmas are easy consequences of the definitions.

\begin{lemma}\label{res:composition_of_interval_preserving_linear_maps}
	The composition of two interval preserving linear maps between vector lattice is again an interval preserving linear map.
\end{lemma}

\begin{lemma}\label{res:interval_preserving_linear_maps_into_sublattices}
	Let $\map\colon E\to F$ be an interval preserving linear map between two vector lattices~$E$ and~$F$. Take a vector sublattice~$F^\prime$ of~$F$ such that $\map(E)\subseteq F^\prime$. Then $\map\colon E\to F^\prime$ is also an interval preserving linear map.
\end{lemma}

The next result, for which we refer to  \cite[Proposition~14.7]{kaplan_THE_BIDUAL_OF_C(X):1988}), is a basic property of interval preserving linear maps.

\begin{proposition}\label{res:interval_preserving_maps_vector_lattices_images}
	Let $\map\colon E\to F$ be an interval preserving linear map between two vector lattices~$E$ and~$F$.
	If~$I$ is an ideal in~$E$, then $\map(I)$ is an ideal in~$F$, and $\pos{\map(I)}=\map(\pos{I})$.
\end{proposition}

We have the following relation between lattice homomorphisms and interval preserving linear maps. The equivalence in~\ref{part:interval_preserving_maps_vector_lattices_1} follows from \cref{res:interval_preserving_maps_vector_lattices_images} and \cite[Proposition~2.1]{van_amstel_van_der_walt_UNPUBLISHED:2022}. The special case of \ref{part:interval_preserving_maps_vector_lattices_2} that an injective interval preserving linear map is a lattice homomorphism was established in \cite[Proposition~2.1]{van_amstel_van_der_walt_UNPUBLISHED:2022}, with a proof different from ours.

\begin{proposition}\label{res:interval_preserving_maps_vector_lattices}
	Let $\map\colon E\to F$ be a linear map between two vector lattices~$E$ and~$F$.
	\begin{enumerate}
		\item\label{part:interval_preserving_maps_vector_lattices_1}
		Suppose that~$\map$ is a lattice homomorphism. Then the following are equivalent:
		\begin{enumerate}
			\item\label{part:interval_preserving_maps_vector_lattices_1_1}
			$\map$ is interval preserving;
			\item\label{part:interval_preserving_maps_vector_lattices_1_2}
			$\map(E)$ is an ideal in~$F$.
		\end{enumerate}	
		\item\label{part:interval_preserving_maps_vector_lattices_2}
		Suppose~$\map$ is interval preserving. Then the following are equivalent:
		\begin{enumerate}
			\item\label{part:interval_preserving_maps_vector_lattices_2_1}
			$\map$ is a lattice homomorphism;
			\item\label{part:interval_preserving_maps_vector_lattices_2_2} $\ker\map$ is an ideal in~$E$.
		\end{enumerate}
	\end{enumerate}
\end{proposition}

\begin{proof}
	In view of the preceding remarks, it only remains to prove that an interval preserving linear map~$\mil$ is a lattice homomorphism when $\ker\mil$ is an ideal in~$F$.  From \cref{res:interval_preserving_maps_vector_lattices_images} we see that $\mil(E)$ is a vector sublattice of~$F$ and that  $\mil(\pos{E})=\pos{\mil(E)}$.
	We let $q\colon E\to E/\ker\mil$ denote the quotient map, and let~$\mil^\prime$ be the linear map $\mil^\prime\colon E/\ker\mil\to F$ such that $\mil=\mil^\prime\circ q$. Then~$\mil^\prime$ is a linear bijection between $E/\ker\mil$ and $\mil(E)$ such that $\mil^\prime(\pos{(E/\ker\mil)})=\mil^\prime(q(\pos{E}))=\mil(\pos{E})=\pos{\mil(E)}$. Hence~$\mil^\prime$ is a lattice homomorphism, and then so is $\mil=\mil^\prime\circ q$.
\end{proof}

The following three results are geared towards direct limits. The verification of the first two is straightforward.

\begin{lemma}\label{res:nested_vector_lattices}
	Let~$E$ be a vector lattice, and let~$I$ be a directed  set. Suppose that, for $i\in I$,~$E_i$ is a vector sublattice of~$E$ such that $E_i\subseteq E_j$ whenever $i\leq j$. Then $\bigcup_i E_i$ is a vector sublattice of~$E$, and the following are equivalent:
	\begin{enumerate}
		\item\label{part:nested_vector_lattices_1}
		all inclusion maps from~$E_i$ into~$E_j$ for $i\leq j$ are interval preserving lattice homomorphisms;
		\item\label{part:nested_vector_lattices_2}
		all inclusion maps from the~$E_i$ into $\bigcup_i E_i$ are interval preserving lattice homomorphisms.
	\end{enumerate}
\end{lemma}

\begin{lemma}\label{res:factoring_map_linear_or_vector_lattice_homomorphism}
	Let~$I$ be a directed  set, let $(E_i)_{i\in I}$ be a collection $\bigl<$~vector spaces / vector lattices~$\bigr>$, and let $\mil_i\colon  E_i\to E$ be $\bigl<$~linear maps / lattice homomorphisms~$\bigr>$ into a $\bigl<$~vector space / vector lattice~$\bigr>$~$E$, such that $\mil_i(E_i)\subseteq\mil_j(E_j)$ when $i\leq j$. Then $\bigcup_i\mil_i(E_i)$ is a $\bigl<$~linear subspace / vector sublattice~$\bigr>$ of~$E$. Let $\factor\colon  \bigcup_i\mil_i(E_i)\to F$ be a map into a $\bigl<$~vector space / vector lattice~$\bigr>$~$F$. Then the following are equivalent:
	\begin{enumerate}
		\item\label{part:factoring_map_vector_lattice_homomorphism_1}
		all maps $\factor\circ\mil_i\colon  E_i\to F$ are $\bigl<$~linear maps / lattice homomorphisms~$\bigr>$;
		\item\label{part:factoring_map_vector_lattice_homomorphism_2}
		$\factor\colon  \bigcup_i\mil_i(E_i)\to F$ is a $\bigl<$~linear map / lattice homomorphism~$\bigr>$.
	\end{enumerate}
\end{lemma}

\begin{proposition}\label{res:factoring_map_interval_preserving}
	Let~$I$ be a directed set, let $(E_i)_{i\in I}$ be a collection vector lattices, and let $\mil_i\colon  E_i\to E$ be interval preserving linear maps into a vector lattice~$E$ such that $\mil_i(E_i)\subseteq\mil_j(E_j)$ when $i\leq j$. Then $\bigcup_i\mil_i(E_i)$ is a vector sublattice of~$E$. Let $\factor\colon  \bigcup_i\mil_i(E_i)\to F$ be a linear map into a vector lattice~$F$. Then the following are equivalent:
	\begin{enumerate}
		\item\label{part:factoring_map_interval_preserving_1}
		all maps $\factor\circ\mil_i\colon  E_i\to F$ are interval preserving linear maps;
		\item\label{part:factoring_map_interval_preserving_2}
		$\factor\colon  \bigcup_i\mil_i(E_i)\to F$ is an interval preserving linear map.	
	\end{enumerate}
\end{proposition}

\begin{proof}
	We know from \cref{res:factoring_map_linear_or_vector_lattice_homomorphism} that  $\bigcup_i\mil_i(E_i)$ is a vector sublattice of~$E$, and that the linearity of all $\factor\circ\mil_i$ is equivalent to that of $\factor$.
	
	We prove that $\factor$ is interval preserving when all $\factor\circ\mil_i$ are. To show that~$\factor$ is positive, take $e\in\pos{\left(\bigcup_i\mil_i(E_i)\right)}$. Then $e=\mil_i(e_i)$ for some~$i$ and $e_i\in E_i$. By \cref{res:interval_preserving_maps_vector_lattices_images},  $\mil_i(e_i)\in\pos{\mil_i(E_i)}=\mil_i(\pos{E_i})$, so we may suppose that $e_i\in\pos{E_i}$. Hence $\factor(e)=(\factor\circ\mil_i)(e_i)\in\pos{F}$, so that~$\factor$ is positive. We show that~$\factor$ is interval preserving. Take $e\in\pos{\left(\bigcup_i\mil_i(E_i)\right)}$ and $f\in[0,\factor(e)]_F$. Again there exist an~$i$ and $e_i\in\pos{E_i}$ such that $\mil_i(e_i)=e$. Since $f\in[0,\factor(e)]_F=[0,(\factor\circ\mil_i)(e_i)]_F$, there exists an $\tilde e_i\in[0,e_i]_{E_i}$ such that $\factor(\mil_i(\tilde e_i))=(\factor\circ\mil_i)(\tilde e_i)=f$. Because $\mil_i(\tilde e_i)\in[0,\mil_i(e_i)]_E=[0,e]_E$, we conclude that~$\factor$ is interval preserving.
	
	It is clear from \cref{res:composition_of_interval_preserving_linear_maps} that all $\factor\circ\mil_i$ are interval preserving when $\factor$ is. 
\end{proof}

The combination of the case of \cref{res:factoring_map_interval_preserving} where~$I$ consists of one element with \cref{res:interval_preserving_maps_vector_lattices_images} and \cref{res:factoring_map_linear_or_vector_lattice_homomorphism} yields the following.

\begin{corollary}\label{res:triangle_lemma_vector_lattices}
	Suppose that the diagram
	\[
	\begin{tikzcd}
		E\arrow{rd}{\maptwo}\arrow{d}[swap]{\map}   &  \\
		\map(E)\arrow[subseteq]{d}\arrow{r}[swap]{\maptwo^\prime}  &   F \\ [-1.1em]
		G&
	\end{tikzcd}
	\]
	is commutative, where~$E$,~$F$, and~$G$ are vector lattices.
	If $\mil\colon E\to G$ and $\maptwo\colon E\to F$ are both $\bigl<$~interval preserving linear maps /  interval preserving lattice homomorphism~$\bigr>$, then $\map(E)$ is a vector sublattice of~$G$, and~$\maptwo^\prime\colon\map(E)\to F$ is an $\bigl<$~interval preserving linear map /  interval preserving lattice homomorphism~$\bigr>$.
\end{corollary}
	
Using \cref{res:composition_of_interval_preserving_linear_maps}, \cref{res:triangle_lemma_vector_lattices} has the following consequence.  We shall apply it in quotient constructions to push down the interval preserving property of the map in the top of a commutative diagram to the map in the bottom. In these applications, $\mil_E$ and $\mil_F$ in the diagram are surjective (quotient) lattice homomorphisms, which are interval preserving according to \cref{res:interval_preserving_maps_vector_lattices}.

\begin{corollary}\label{res:square_lemma_vector_lattices}
	Suppose that the diagram
	\[
	\begin{tikzcd}
		E\arrow{r}{\maptwo}\arrow[two heads]{d}[swap]{\map_E}   &  F\arrow{d}{\map_F}\\
		E^\prime\arrow{r}[swap]{\maptwo^\prime}  &   F^\prime
	\end{tikzcd}
	\]	
	is commutative, where~$E$,~$E^\prime$,~$F$, and~$F^\prime$ are vector lattices;~and~$\map_E$ is surjective.
	If~$\mil_E$,~$\mil_F$, and~$\maptwo$ are $\bigl<$~interval preserving linear maps / interval preserving lattice homomorphism~$\bigr>$, then~$\maptwo^\prime$ is an $\bigl<$~interval preserving linear map / interval preserving lattice homomorphism~$\bigr>$.
\end{corollary}


\subsection{Almost interval preserving linear maps and lattice homomorphisms}\label{subsec:almost_interval_preserving_linear_maps}


The connection between (almost) interval preserving linear maps and lattice homomorphisms is particularly strong in the case of normed lattices and continuous maps.
We recall that the norm dual~$\ndual{E}$  of a normed lattice~$E$ is an ideal in~$\odual{E}$ that is a Banach lattice, and that $\ndual{E}=\odual{E}$ when~$E$ is a Banach lattice; see \cite[Theorem~3.49, Theorem~4.1, and Corollary~4.5]{aliprantis_burkinshaw_POSITIVE_OPERATORS_SPRINGER_REPRINT:2006}, for example. We then have the following, which is for the most part \cite[Proposition~1.4.19]{meyer-nieberg_BANACH_LATTICES:1991}. Part~\ref{part:duality_2_3}, which applies to order continuous Banach lattices, has been added for use in the present paper. It also brings out the perfect symmetry between~$\mil$ and~$\ndual{\mil}$ when order intervals are weakly compact.

\begin{proposition}\label{res:duality}
	Let $\mil\colon  E\to F$ be a continuous linear map between normed lattices~$E$ and~$F$, with adjoint $\ndual\mil\colon  \ndual{F}\to\ndual{E}$.
	\begin{enumerate}
		\item\label{part:duality_1}
		The following are equivalent:
		\begin{enumerate}
			\item\label{part:duality_1_1}
			$\mil\colon  E\to F$ is a lattice homomorphism;
			\item\label{part:duality_1_2}
			$\ndual{\mil}\colon \ndual{F}\to\ndual{E}$ is almost interval preserving;
			\item\label{part:duality_1_3}
			$\ndual{\mil}\colon  \ndual{F}\to\ndual{E}$ is interval preserving.
		\end{enumerate}
		\item\label{part:duality_2}
		The following are equivalent:
		\begin{enumerate}
			\item\label{part:duality_2_1}
			$\ndual{\mil}\colon  \ndual{F}\to\ndual{E}$ is a lattice homomorphism;
			\item\label{part:duality_2_2}
			$\mil\colon  E\to F$ is almost interval preserving.
		\end{enumerate}
		When order intervals in~$E$ are weakly compact, then these are also equivalent to:
		\begin{enumerate}[resume]
			\item\label{part:duality_2_3}
			$\mil\colon  E\to F$ is interval preserving.
		\end{enumerate}
	\end{enumerate}
\end{proposition}
\begin{proof}
	In view of the results in \cite[Proposition~1.4.19]{meyer-nieberg_BANACH_LATTICES:1991}, we need only prove that~\ref{part:duality_2_2} implies~\ref{part:duality_2_3} when order intervals in~$E$ are weakly compact. Take $x\in\pos{E}$. Since the order interval $[0,x]$ is weakly compact, so is $\mil([0,x])$. Hence it is weakly closed. Because it is convex, this implies that $\mil([0,x])$ is norm closed. Since~$\mil$ is supposed to be almost interval preserving, we now see that it is actually even interval preserving.
\end{proof}

The following two lemmas are direct consequences of the definitions.

\begin{lemma}\label{res:composition_of_almost_interval_preserving_linear_maps}
	Let~$E$ be a vector lattice, and let~$F$ and~$G$ be normed lattices.  Suppose that $\mil\colon  E\to F$ is an almost interval preserving linear map, and that $\maptwo\colon  F\to G$ is a continuous almost interval preserving linear map. Then $\maptwo\circ\mil\colon  E\to G$ is an almost interval preserving linear map.
\end{lemma}

\begin{lemma}\label{res:almost_interval_preserving_linear_maps_into_sublattices}
	Let $\mil\colon  E\to F$ be an almost interval preserving linear map between a vector lattice~$E$ and a normed lattice~$F$. Take a vector sublattice~$F^\prime$ of~$F$ such that $\mil(E)\subseteq F^\prime$. Then $\mil\colon  E\to F^\prime$ is also an almost interval preserving linear map.
\end{lemma}

It is established in \cite[Theorem~2.1.(1)]{bouras_elbour:2018} that $\overline{\mil(E)}$ is an ideal in~$F$ when $\mil\colon  E\to F$ is an almost interval preserving linear map between Banach lattices~$E$ and~$F$. Our next result is more precise. For its proof, we recall from \cref{sec:preliminaries} that, for a subset $S$ of a vector lattice $E$, $\pos{S}$ denotes the set of positive \emph{parts} of elements of $S$. This contains $S\cap \pos{E}$, but the inclusion can be proper.

\begin{proposition}\label{res:almost_interval_preserving_maps_normed_vector_lattices}
	Let $\map\colon E\to F$ be an almost interval preserving linear map between a vector lattice~$E$ and a normed lattice~$F$. Let~$I$ be an ideal in~$E$. Then:
	\begin{enumerate}
		\item\label{part:almost_interval_preserving_maps_normed_vector_lattices_1}
		$\overline{\map(I)}$ is an ideal in~$F$ and $\pos{\left(\overline{\map(I)}\right)}=\overline{\map(\pos{I})}=\overline{\pos{\mil(I)}}$;
		\item\label{part:almost_interval_preserving_maps_normed_vector_lattices_2}
		the following are equivalent:
		\begin{enumerate}
			\item\label{part:almost_interval_preserving_maps_normed_vector_lattices_2_1}
 			$\mil(\pos{I})$ is closed;
 			\item\label{part:almost_interval_preserving_maps_normed_vector_lattices_2_2}
 			$\mil(I)$ is closed and $\pos{\mil(I)}=\map(\pos{I})$.
 		\end{enumerate}
 	\end{enumerate}
\end{proposition}

\begin{proof}	
	We prove~\ref{part:almost_interval_preserving_maps_normed_vector_lattices_1}. As a preparation, we establish the following claim: whenever $y,z\in F$ are such that $0\leq y\leq \abs{z}$ and $z\in\overline{\map(I)}$, then $y\in \overline{\map(\pos{I})}$. To see this, choose a sequence $(x_n)\subseteq I$ such that $\map(x_n)\to z$. Then $\abs{\map(x_n)}\to \abs{z}$. Since $\abs{\map(x_n)}\in[0,\map(\abs{x_n})]_F=\overline{\map([0,\abs{x_n}]_E)}$, there exists a sequence $(x_n^\prime)$ in~$E$ such that $x_n^\prime\in[0,\abs{x_n}]_E$ and $\norm{\map(x_n^\prime)-\abs{\map (x_n)}}<1/2^n$. Then $(x_n^\prime)\subseteq\pos{I}$ and $\map(x_n^\prime)\to\abs{z}$. Hence $y\wedge \map(x_n^\prime)\to y\wedge\abs{z}=y$. Since $y\wedge\map(x_n^\prime)\in[0,\map(x_n^\prime)]_F=\overline{\map([0,x_n^\prime]_E)}$, there exists a sequence $(x_n^{\prime\prime})$ in~$E$ such that $x_n^{\prime\prime}\in[0,x_n^\prime]_E$ and $\norm{\map(x_n^{\prime\prime})-y\wedge \map(x_n^\prime)}<1/2^n$. Then $(x_n^{\prime\prime})\subseteq \pos{I}$ and $\map(x_n^{\prime\prime})\to y$, showing that $y\in\overline{\map(\pos{I})}$. Our claim has now been established.
	
	We can now show that $\overline{\map(I)}$ is an ideal in~$F$. Suppose that $y,z\in F$ are such that $z\in\overline{\map(I)}$ and that $0\leq\abs{y}\leq\abs{z}$. Then $0\leq y^\pm\leq\abs{z}$. It follows from the claim that $y^\pm\in \overline{\map(\pos{I})}\subseteq\overline{\map(I)}$, so that also $y\in\overline{\map(I)}$. Hence $\overline{\map(I)}$ is an ideal in~$F$. 
	
	We show that $\pos{\left(\overline{\map(I)}\right)}=\overline{\map(\pos{I})}$. Take $z\in\pos{\left(\overline{\map(I)}\right)}$. Since we already know that $\overline{\map(I)}$ is a vector sublattice of~$F$, we have $z\in\overline{\map(I)}$. From the fact that $0\leq z\leq \abs{z}$ and the claim it follows that $z\in\overline{\map(\pos{I})}$. Hence $\pos{\left(\overline{\map(I)}\right)}\subseteq\overline{\map(\pos{I})}$. The converse inclusion $\overline{\map(\pos{I})}\subseteq\pos{\left(\overline{\map(I)}\right)}$ follows from the positivity of $\map$ and the closedness of~$\pos{F}$.
	
	Next we show that $\overline{\map(\pos{I})}=\overline{\pos{\mil(I)}}$. The positivity of~$\mil$ makes clear that $\overline{\map(\pos{I})}\subseteq \overline{\pos{\mil(I)}}$. For the reverse inclusion, take $z\in \overline{\pos{\mil(I)}}$. There exists a sequence $(x_n)$ in~$I$ such that $\pos{\mil(x_n)}\to z$. Since $0\leq\pos{\mil(x_n)}\leq\abs{\mil(x_n)}\leq\mil(\abs{x_n})$, we have $\pos{\mil(x_n)}\in[0,\mil(\abs{x_n})]_F$. From this we see that there exist $x_n^\prime\in[0,\abs{x_n}]_E\subseteq \pos{I}$ such that $\norm{\mil(x_n^\prime)-\pos{\mil(x_n)}}<1/2^n$. Then also $\mil(x_n^\prime)\to z$, showing that $z\in\overline{\map(\pos{I})}$. The proof of~\ref{part:almost_interval_preserving_maps_normed_vector_lattices_1} is now complete.
	
	We prove that~\ref{part:almost_interval_preserving_maps_normed_vector_lattices_2_1} implies~\ref{part:almost_interval_preserving_maps_normed_vector_lattices_2_2}. In view of~\ref{part:almost_interval_preserving_maps_normed_vector_lattices_1}, it is sufficient to prove that $\mil(I)$ is closed. For this, suppose that $z\in F$ and that $(x_n)\subseteq I$ is a sequence such that $\mil(x_n)\to z$. Then $\mil(x_n)^\pm\to z^\pm$, so that $z^\pm\in\overline{\pos{\mil(I)}}$. By~\ref{part:almost_interval_preserving_maps_normed_vector_lattices_1}, the latter set equals $\overline{\map(\pos{I})}$. Since $\map(\pos{I})$ is closed, we see that $z^\pm\in\map(\pos{I})\subseteq\map(I)$. Then also $z\in\map(I)$, as desired.
	
	It is immediate from~\ref{part:almost_interval_preserving_maps_normed_vector_lattices_1} that~\ref{part:almost_interval_preserving_maps_normed_vector_lattices_2_2} implies~\ref{part:almost_interval_preserving_maps_normed_vector_lattices_2_1}.
\end{proof}

The following result is concerned with the relation between lattice homomorphisms and almost interval preserving linear maps.

\begin{proposition}\label{res:almost_interval_preserving_maps_and_lattice_homomorphisms_normed_vector_lattices}
Let $\map\colon E\to F$ be a linear map between a vector lattice~$E$ and a normed lattice~$F$.
	\begin{enumerate}
		\item\label{part:almost_interval_preserving_maps_and_lattice_homomorphisms_normed_vector_lattices_1}
		Suppose that~$\map$ is a lattice homomorphism. Then the following are equivalent:
			\begin{enumerate}
				\item\label{part:almost_interval_preserving_maps_and_lattice_homomorphisms_normed_vector_lattices_1_1}
				$\map$ is almost interval preserving;
				\item\label{part:almost_interval_preserving_maps_and_lattice_homomorphisms_normed_vector_lattices_1_2}
				$\overline{\map(E)}$ is an ideal in~$F$.
			\end{enumerate}
		\item\label{part:almost_interval_preserving_maps_and_lattice_homomorphisms_normed_vector_lattices_2}
		Suppose that~$\mil$ is almost interval preserving.
		If, in addition, $\ker \mil$ is an ideal in~$E$ and $\mil(\pos{E})$ is closed, then~$\mil$ is a lattice homomorphism.
	\end{enumerate}
\end{proposition}

\begin{proof}
	We prove~\ref{part:almost_interval_preserving_maps_and_lattice_homomorphisms_normed_vector_lattices_1}. In view of \cref{res:almost_interval_preserving_maps_normed_vector_lattices}, it only remains to be shown that~$\map$ is almost interval preserving when~$\map$ is a lattice homomorphism and $\overline{\map(E)}$ is an ideal in~$F$. In this case, take~$x\in\pos{E}$ and $y\in F$ such that $y\in [0,\map(x)]_F$. Since $\map(x)\in\overline{\map(E)}$, we also have $y\in \overline{\map(E)}$. Choose a sequence $(x_n)$ in~$E$ such that $\map(x_n)\to y$. Then $y=|y|\land \map(x)=\left[\lim_{n\to\infty}|\map(x_n)|\right]\land \map(x)=\lim_{n\to\infty}\map(|x_n|\land x)$. Since $|x_n|\land x\in [0,x]$, this shows that $y\in\overline{\map([0,x]_E)}$.

	We prove~\ref{part:almost_interval_preserving_maps_and_lattice_homomorphisms_normed_vector_lattices_2}.
		We see from \cref{res:almost_interval_preserving_maps_normed_vector_lattices} that $\mil(E)$ is a vector sublattice of~$F$ and that $\mil(\pos{E})=\pos{\mil(E)}$. The argument in the proof of part~\ref{part:interval_preserving_maps_vector_lattices_2} of \cref{res:interval_preserving_maps_vector_lattices} then shows that~$\mil$ is a lattice homomorphism.
	\end{proof}

In particular, the inclusion map from a vector sublattice of a normed lattice into its closure is almost interval preserving.

If~$E$ is a normed lattice where order intervals are weakly compact and if~$\mil$ is continuous, then part~\ref{part:almost_interval_preserving_maps_and_lattice_homomorphisms_normed_vector_lattices_2} of \cref{res:almost_interval_preserving_maps_and_lattice_homomorphisms_normed_vector_lattices} can be improved.

\begin{proposition}\label{res:almost_interval_preserving_linear_maps_and_homomorphism_order_intervals_weakly_compact}
	Let~$\mil\colon E\to F$ be a continuous almost interval preserving linear map between normed vector lattices~$E$ and~$F$, where order intervals in~$E$ are weakly compact. Then the following are equivalent:
	\begin{enumerate}
		\item\label{part:almost_interval_preserving_linear_maps_and_homomorphism_order_intervals_weakly_compact_1}
		$\mil$ is a lattice homomorphism;
		\item\label{part:almost_interval_preserving_linear_maps_and_homomorphism_order_intervals_weakly_compact_2}
		$\ker \mil$ is an ideal in~$E$.
	\end{enumerate}
\end{proposition}

\begin{proof}
	We need only prove that~\ref{part:almost_interval_preserving_linear_maps_and_homomorphism_order_intervals_weakly_compact_2} implies~\ref{part:almost_interval_preserving_linear_maps_and_homomorphism_order_intervals_weakly_compact_1}.
	By \cref{res:duality},~$\mil$ is, in fact, even interval preserving, and then part~\ref{part:interval_preserving_maps_vector_lattices_2} of \cref{res:interval_preserving_maps_vector_lattices} shows that~$\mil$ is a lattice homomorphism.
\end{proof}	

The majority of the results in the remainder of this section are relevant in the context of direct limits in $\NLAIPLH$ and $\BLAIPLH$.
In view of part~\ref{part:almost_interval_preserving_maps_normed_vector_lattices_1} of \cref{res:almost_interval_preserving_maps_normed_vector_lattices}, the first one is about closures of vector sublattices being ideals.

\begin{proposition}\label{res:nested_normed_vector_lattices}
	Let~$E$ be a normed lattice, and let~$I$ be a directed set. Suppose that, for $i\in I$,~$E_i$ is a vector sublattice of~$E$ such that $E_i\subseteq E_j$ whenever $i\leq j$, and that $E=\overline{\bigcup_i E_i}$. Then the following are equivalent:
	\begin{enumerate}
		\item\label{part:nested_normed_vector_lattices_1}
		all inclusion maps from~$E_i$ into~$E_j$ for $i\leq j$ are almost interval preserving lattice homomorphisms;
		\item\label{part:nested_normed_vector_lattices_2}
		all inclusion maps from the~$E_i$ into the vector sublattice $\bigcup_i E_i$ of~$E$ are almost interval preserving lattice homomorphisms;
		\item\label{part:nested_normed_vector_lattices_3}
		all inclusion maps from the~$E_i$ into~$E$ are almost interval preserving lattice homomorphisms.
	\end{enumerate}
\end{proposition}

\begin{proof}
	We prove that~\ref{part:nested_normed_vector_lattices_1} implies~\ref{part:nested_normed_vector_lattices_3}. Fix an index~$i$ and an $e_i\in\pos{E}$, and take $e\in[0,e_i]_{E}$.  There exist a sequence of indices $(i_n)$ and elements $e_{i_n}$ of $E_{i_n}$ such that $e_{i_n}\to e$. We may suppose that $i_n\geq i$  and that $e_{i_n}\in\pos{E_{i_n}}$ for all~$n$. We also have that $e_{i_n}\wedge e_i\to e$. Since $e_{i_n}\wedge e_i\in[0,e_i]_{E_{i_n}}$, there exists a sequence $(e_{i_n}^\prime)$ such that $e_{i_n}^\prime\in[0,e_i]_{E_i}$ and $\norm{e_{i_n}^\prime-e_{i_n}\wedge e_i}< 1/2^n$ for all~$n$. This implies that also $e_{i_n}^\prime\to e$. Hence $e\in\overline{[0,e_i]}^{E}$, as desired.
	
	It is immediate from \cref{res:almost_interval_preserving_linear_maps_into_sublattices} that
	\ref{part:nested_normed_vector_lattices_3} implies~\ref{part:nested_normed_vector_lattices_2} and that~\ref{part:nested_normed_vector_lattices_2} implies~\ref{part:nested_normed_vector_lattices_1}.
\end{proof}

\newpage 

The proof of the following result is elementary.

\begin{lemma}\label{res:factoring_map_linear_or_vector_lattice_homomorphism_normed_case}
	Let~$I$ be a directed set, let $(E_i)_{i\in I}$ be a collection $\bigl<$~vector spaces / vector lattices~$\bigr>$, and let $\mil_i\colon  E_i\to E$ be $\bigl<$~linear maps / lattice homomorphisms~$\bigr>$ into a $\bigl<$~normed space / normed lattice~$\bigr>$~$E$, such that $\mil_i(E_i)\subseteq\mil_j(E_j)$ when $i\leq j$. Let $\factor\colon \overline{\bigcup_i\mil_i(E_i)}  \to F$ be a continuous map into a $\bigl<$~normed space / normed lattice~$\bigr>$~$F$. Then $\overline{\bigcup_i\mil_i(E_i)}$ is a $\bigl<$~vector space / vector lattice~$\bigr>$, and the following are equivalent:
	\begin{enumerate}
		\item\label{part:factoring_map_vector_lattice_homomorphism_normed_case_1}
		all maps $\factor\circ\mil_i\colon  E_i\to F$ are $\bigl<$~linear maps / lattice homomorphisms~$\bigr>$;
		\item\label{part:factoring_map_vector_lattice_homomorphism_normed_case_3}
		$\factor$ is a $\bigl<$~linear map / lattice homomorphism~$\bigr>$.
	\end{enumerate}
\end{lemma}

The proof of our next result is less straightforward than that of its counterpart for interval preserving linear maps between vector lattices, \cref{res:factoring_map_interval_preserving}.

\begin{proposition}\label{res:factoring_map_almost_interval_preserving}
	Let~$I$ be a directed set, let $(E_i)_{i\in I}$ be a collection of vector lattices, let $\mil_i\colon  E_i\to E $ be almost interval preserving linear maps into a normed lattice~$E$ such that $\mil_i(E_i)\subseteq\mil_j(E_j)$. Then $\overline{\bigcup_i\mil_i(E_i)}$ is a vector sublattice of~$E$.  Let $\factor\colon \overline{\bigcup_i\mil_i(E_i)}\to F$ be a continuous linear map into a normed lattice~$F$. Then the following are equivalent:
	\begin{enumerate}
		\item\label{part:factoring_map_almost_interval_preserving_1}
		all maps $\factor\circ\mil_i\colon  E_i\to F$ are almost interval preserving linear maps;
		\item\label{part:factoring_map_almost_interval_preserving_2}
		$\factor$ is an almost interval preserving linear map.	
	\end{enumerate}
\end{proposition}

\begin{proof}
	We know from \cref{res:factoring_map_linear_or_vector_lattice_homomorphism_normed_case} that  $\overline{\bigcup_i\mil_i(E_i)}$ is a vector sublattice of~$E$, and that the linearity of all $\factor\circ\mil_i$ is equivalent to that of $\factor$.
	
	We prove that $\factor$ is almost interval preserving when all $\factor\circ\mil_i$ are.
	We prove that~\ref{part:factoring_map_almost_interval_preserving_1} implies~\ref{part:factoring_map_almost_interval_preserving_2}. To show that~$\factor$ is positive, let $e\in\pos{\Big(\overline{\bigcup_i\mil_i(E_i)}\Big)}$. There exist sequences $(i_n)$ of indices and $(e_{i_n})$ of elements $e_{i_n}$ of $E_{i_n}$ such that $\mil_{i_n}(e_{i_n})\to e$. Then also $(\mil_{i_n}(e_{i_n}))^+\to e$. As in the proof of part~\ref{part:almost_interval_preserving_maps_normed_vector_lattices_1} of \cref{res:almost_interval_preserving_maps_normed_vector_lattices}, we have $(\mil_{i_n}(e_{i_n}))^+\in[0,\mil(\abs{e_{i_n}})]_F$, so that there exist $\tilde e_{i_n}\in \pos{(E_{i_n})}$ such that $\norm{\mil_{i_n}(\tilde e_{i_n})-(\mil_{i_n}(e_{i_n}))^+}<1/2^n$. Since then also $\mil_{i_n}(\tilde e_{i_n})\to e$, we see from $\factor(e)=\lim_n(\factor\circ\mil_{i_n})(\tilde e_{i_n})$ that~$\factor$ is positive.
	
	To see that~$\factor$ is almost interval preserving when the $\factor\circ \mil_i$ are, take $e\in\pos{\Big(\overline{\bigcup_i\mil_i(E_i)}\Big)}$. Suppose that $f\in [0,\factor(e)]_F$. Take $\veps>0$. There exist an index~$i$ and an $e_i\in E_i$ such that $\lrnorm{e-\mil_i(e_i)}<\veps$. Using that then also $\lrnorm{e-\pos{\mil_i(e_i)}}<\veps$, we see from the second equality in part~\ref{part:almost_interval_preserving_maps_normed_vector_lattices_1} of  \cref{res:almost_interval_preserving_maps_normed_vector_lattices} that we may suppose that $e_i\in\pos{E_i}$. Set $\widetilde f\coloneqq f\wedge (\factor\circ \mil_i) (e_i)\in [0,(\factor\circ\mil_i) (e_i)]_F$.  Hence there exists an $\widetilde e_i\in [0,e_i]_{E_i}$ such that $\big\lVert\widetilde f-(\factor\circ\mil_i)(\widetilde e_i)\big\rVert<\veps$, and then we have
	\begin{align*}
		\big\lVert f-\factor\big(\mil_i (\widetilde e_i)\wedge e\big)\big\rVert & \leq \big\lVert f-\widetilde f\big\rVert +\big\lVert\widetilde f-(\factor\circ\mil_i) (\widetilde e_i)\big\rVert\\
		&\phantom{=}\quad+\big\lVert(\factor\circ\mil_i) (\widetilde e_i)-\factor\big(\mil_i (\widetilde e_i)\wedge e\big)\big\rVert\\
		&= \lrnorm{f\wedge \factor(e)-f\wedge\factor\big(\mil_i (e_i)\big)}+\big\lVert\widetilde f-(\factor\circ\mil_i) (\widetilde e_i)\big\rVert\\
		&\phantom{=}\quad+\big\lVert\factor\big(\mil_i (\widetilde e_i)\wedge\mil_i (e_i)\big)-\factor\big(\mil_i( \widetilde e_i)\wedge e\big)\big\rVert\\
		&< \big\lVert\factor(e)-\factor\big(\mil_i (e_i)\big)\big\rVert+\veps\\
		&\phantom{=}\quad+\big\lVert\factor\big\rVert \, \big\lVert \mil_i (\widetilde e_i)\wedge\mil_i (e_i)-\mil_i (\widetilde e_i)\wedge e\big\rVert\\
		&\leq \lrnorm{\factor} \lrnorm{e-\mil_i (e_i)}+\veps+\lrnorm{\factor}\,\lrnorm{\mil_i (e_i)- e}\\
		&\leq (1+2\lrnorm{\factor})\,\veps.
	\end{align*}	
	As $\mil_i (\widetilde e_i)\wedge e\in [0,e]_{\overline{\bigcup_i\mil_i(E_i)}}$, we conclude that $f\in\overline{\factor\big([0,e]_{\overline{\bigcup_i\mil_i(E_i)}}\big)}$, as desired.
	
	\cref{res:composition_of_almost_interval_preserving_linear_maps} and \cref{res:almost_interval_preserving_linear_maps_into_sublattices} show that~all $\factor\circ\mil_i$ are almost interval preserving when $\factor$ is.
\end{proof}

When the~$E_i$ in \cref{res:factoring_map_almost_interval_preserving} are normed lattices and the~$\mil_i$ are continuous, \cref{res:duality} provides the ingredients for an alternate proof of a variation of \cref{res:factoring_map_almost_interval_preserving}. Although the result is weaker, we still include it for reasons of aesthetic appeal of its proof.

\begin{proposition}\label{res:factoring_map_almost_interval_preserving_not_nested}
	Let~$I$ be a set, let $(E_i)_{i\in I}$ be a collection normed lattices, and let $\mil_i\colon E_i\to E$ be continuous almost interval preserving linear maps into a normed lattice~$E$, such that $E=\overline{\bigcup_i\mil_i(E_i)}$. Let $\factor\colon E\to F$ be a continuous linear map into a normed lattice~$F$. Then the following are equivalent:
	\begin{enumerate}
		\item\label{part:factoring_map_almost_interval_preserving_not_nested_1}
		all maps $\factor\circ\mil_i\colon E_i\to F$ are almost interval preserving linear maps;
		\item\label{part:factoring_map_almost_interval_preserving_not_nested_2}
		$\factor$ is an almost interval preserving linear map.	
	\end{enumerate}
\end{proposition}

\begin{proof}
	To show that $\factor$ is almost interval preserving when the $\factor\circ\mil_i$ are, we use the equivalence of~\ref{part:duality_2_1} and~\ref{part:duality_2_2} of \cref{res:duality}. We know that the $\ndual{\mil_i}\colon \ndual{E}\to\ndual{E_i}$ and the $\ndual{(\factor\circ\mil_i)}\colon\ndual{F}\to\ndual{E_i}$ are all lattice homomorphisms and need to show that $\ndual{\factor}\colon\ndual{F}\to\ndual{E}$ is also a lattice homomorphism. Take $\ndual{f}\in\ndual{F}$. Then
	\begin{align*}
		\ndual{\mil_i}\big[\ndual{\factor}(\abs{\ndual f})\big]&=\ndual{(\factor\circ\mil_i)}(\abs{\ndual{f}})=\lrabs{\ndual{(\factor\circ\mil_i)}(\ndual{f})}\\
		&=\lrabs{\ndual{\mil_i}\big[\ndual{\factor}(\ndual{f})\big]}=\ndual{\mil_i}\big[\lrabs{\ndual{\factor}(\ndual{f})}\big].
	\end{align*}
	Stated otherwise, we have that $\big[\ndual{\factor}(\abs{\ndual{f}})\big](\mil_i(e_i))=\big[\lrabs{\ndual{\factor}(\ndual{f})}\big](\mil_i(e_i))$
	for $e_i\in E_i$. Since $\bigcup_i\mil_i(E_i)$ is dense in~$E$, we conclude that  $\ndual{\factor}(\abs{\ndual{f}})=\lrabs{\ndual{\factor}(\ndual{f})}$. Hence~$\ndual{\factor}$ is a lattice homomorphism.
	
		\cref{res:composition_of_almost_interval_preserving_linear_maps} shows that~all $\factor\circ\mil_i$ are almost interval preserving when $\factor$ is.
\end{proof}
	
We have already observed that the inclusion map from a vector sublattice of a normed lattice into its closure is almost interval preserving. Hence the case of \cref{res:factoring_map_almost_interval_preserving} where~$I$ has one element has the following consequence.

\begin{corollary}\label{res:factoring_map_almost_interval_preserving_dense_sublattice}	
	Let~$E^\prime$ be a dense vector sublattice of a normed lattice~$E$, let~$F$ be a normed lattice, and let $\factor\colon E\to F$ be a continuous linear map. Then~$\factor$ is almost interval preserving if and only if its restriction $\factor\colon E^\prime\to F$ to~$E^\prime$ is almost interval preserving.
\end{corollary}	
	
The combination of the case of \cref{res:factoring_map_almost_interval_preserving} where~$I$ consists of one element with part~\ref{part:almost_interval_preserving_maps_normed_vector_lattices_1} of \cref{res:almost_interval_preserving_maps_normed_vector_lattices} and with  \cref{res:factoring_map_linear_or_vector_lattice_homomorphism_normed_case} yields the following companion result to \cref{res:triangle_lemma_vector_lattices}.

\begin{corollary}\label{res:triangle_lemma_normed_vector_lattices}
	Suppose that the diagram
	\[
	\begin{tikzcd}
		E\arrow{rd}{\maptwo}\arrow{d}[swap]{\map}   &  \\
		\overline{\map(E)}\arrow[subseteq]{d}\arrow{r}[swap]{\maptwo^\prime}  &   F \\ [-1.1em]
		G&
	\end{tikzcd}
	\]
	is commutative, where~$E$ is vector lattice;~$F$ and~$G$ are normed lattices; and~$\maptwo^\prime$ is continuous. If $\mil\colon E\to G$ and $\maptwo\colon E\to F$ are both $\bigl<$~almost interval preserving linear maps /  almost interval preserving lattice homomorphism~$\bigr>$, then $\overline{\map(E)}$ is a vector sublattice of~$G$, and~$\maptwo^\prime\colon\overline{\map(E)}\to F$ is an $\bigl<$~almost interval preserving linear map /  almost interval preserving lattice homomorphism~$\bigr>$.
\end{corollary}

Using \cref{res:composition_of_almost_interval_preserving_linear_maps}, \cref{res:triangle_lemma_normed_vector_lattices} is easily seen to have the following consequence, which is a companion result to \cref{res:square_lemma_vector_lattices}. Analogously to the latter result, we shall apply it in quotient constructions to push down the almost interval preserving property from the map in the top of a commutative diagram to the map in the bottom, in cases where~$\mil_E$ and~$\mil_F$ in the diagram are surjective quotient lattice homomorphisms.

\begin{corollary}\label{res:square_lemma_normed_vector_lattices}
	Suppose that the diagram
	\[
	\begin{tikzcd}
		E\arrow{r}{\maptwo}\arrow{d}[swap]{\mil_E}   &  F\arrow{d}{\mil_F}\\
		E^\prime\arrow{r}[swap]{\maptwo^\prime}  &   F^\prime
	\end{tikzcd}
	\]	
	is commutative, where~$E$ is a vector lattice;~$E^\prime$,~$F$ and~$F^\prime$ are normed lattices; $\overline{\mil_E(E)}=E^\prime$; and~$\mil_F$ and~$\maptwo^\prime$ are continuous.
	If~$\mil_E$,~$\mil_F$, and~$\maptwo$ are $\bigl<$~almost interval preserving linear maps / almost interval preserving lattice homomorphism~$\bigr>$, then~$\maptwo^\prime$ is an $\bigl<$~almost interval preserving linear map / almost interval preserving lattice homomorphism~$\bigr>$.
\end{corollary}


\section{Direct limits: existence via three standard constructions}\label{sec:direct_limits_existence_via_three_standard_constructions}


\noindent In each of the categories $\VS$, $\NS$, and $\BS$, every direct system has a direct limit. This can be proved via well-known standard constructions that are similar to each other. We shall now show that each of these constructions can also be used to provide direct limits in some (but not all) of the categories of vector lattices, normed latticed, and Banach lattices under consideration in this paper.


\subsection{Direct limits of vector lattices}\label{subsec:direct_limits_of_vector_lattices}


Every direct system in the category $\VLLH$ of vector lattices and lattice homomorphisms has a direct limit in $\VLLH$. This was first established by Filter in \cite{filter:1988}. His method is to first regard the system as a direct system in $\SET$, and then supply the canonical direct limit in this category with the structure of a vector lattice to produce a direct limit in $\VLLH$. This is also the method as reviewed in \cite{van_amstel_van_der_walt_UNPUBLISHED:2022}. There is, however, an alternate approach that is, perhaps, more transparent. It uses a well-known construction that is applied in many non-analytical categories of practical interest. The idea is to start with a direct system in the category $\VS$ of vector spaces and linear maps, and construct a direct limit of that system in $\SET$, but in such a way that the set~$E$ in that direct limit in $\SET$ is already naturally a vector space. It is then virtually immediately clear that also a direct limit of the system in $\VS$ has been found. If the system is in $\VLLH$, then this standard construction (as we shall call it) equally evidently produces a direct limit in $\VLLH$. With a little help from the general toolbox in \cref{sec:interval_preserving_and_almost_interval_preserving_linear_maps_and_lattice_homomorphisms}, this also works in the category $\VLIPLH$ of vector lattices and interval preserving lattice homomorphisms, where the existence of direct limits was first established in \cite[Proposition~3.4]{van_amstel_van_der_walt_UNPUBLISHED:2022}).

The whole setup shows naturally that certain direct limits are also direct limits in categories `to the right' in a chain of categories, as in \cite[Proposition~3.4]{van_amstel_van_der_walt_UNPUBLISHED:2022}. Furthermore, in contrast to the approach in \cite{filter:1988}, this construction can easily be adapted to work for direct limits of normed lattices and Banach lattices.

Suppose, then, that $\Edirsys$ is a direct system in $\VS$. We construct a direct limit of the system in $\SET$. For this, consider the vector space
\[
	{\widetilde E}\coloneqq\prod_i E_i
\]
and its linear subspace
\[
	{\widetilde E}_0\coloneqq\left\{(e_i)\in {\widetilde E}: \text{there exists an index } i\text{ such that } e_j=0 \text{ for }j\geq i\right\}.
\]
We define the linear maps $\maptwo_i\colon E_i\to {\widetilde E}$ by setting
\begin{equation}\label{eq:psi_map_definition}
	\bigl(\maptwo_i(e_i)\big)_j\coloneqq
	\begin{cases}
		\cm_{ji}(e_i)&\text{when }j\geq i;\\
		0&\text{else}.
	\end{cases}
\end{equation}
We let $q\colon {\widetilde E}\to {\widetilde E}/{\widetilde E}_0$ denote the linear quotient map into the vector space $\widetilde E/\widetilde E_0$, and define the linear maps $\mil_i\colon E_i\to{\widetilde E}/{\widetilde E}_0$ by setting $\mil_i\coloneqq q\circ\maptwo_i$. Clearly, when $e_i\in E_i$ and $e_j\in E_j$, then $\mil_i(e_i)=\mil_j(e_j)$ if and only if there exists a $k\geq i,j$ such that $\cm_{ki}(e_i)=\cm_{kj}(e_j)$. It is then immediate that $\mil_i=\mil_j\circ \cm_{ji}$ when $i\leq j$. We define the subset
\[
	E\coloneqq\bigcup_i\mil_i(E_i)
\]
of the space ${\widetilde E}/{\widetilde E}_0$. The fact that $\mil_i=\mil_j\circ\cm_{ji}$ for $i\leq j$ implies that $\mil_i(E_i)\subseteq\mil_j(E_j)$ when $j\geq i$. Hence~$E$ is a nested union of linear subspaces of ${\widetilde E}/{\widetilde E}_0$, so that~$E$ itself is also a linear subspace of  ${\widetilde E}/{\widetilde E}_0$; we view the~$\mil_i$ as linear maps from the~$E_i$ into~$E$.

We claim that $\Esys$ is a direct limit of $\Edirsys$ in $\SET$. The compatibility with $\Edirsys$ has already been observed. For the universal property, suppose that $\Esysprimed$ is a system in $\SET$ that is compatible with $\Edirsys$. If $\factor\colon E\to E^\prime$ is to be a factoring map, then the requirement $\mil_i^\prime=\factor\circ\mil_i$ obviously uniquely determines~$\factor$. Hence we take this as a definition: for $e\in E$, choose~$i$ and $e_i\in E_i$ such that $e=\mil_i(e_i)$, and set $\factor(e)\coloneqq\mil_i^\prime(e_i)$. To show that this is well defined, suppose that $\mil_i(e_i)=\mil_j(e_j)$. Then there exists a $k\geq i,j$ such that $\cm_{ki}(e_i)=\cm_{kj}(e_j)$, which implies that $\mil_i^\prime(e_i)=\mil_k^\prime\bigl(\cm_{ki}(e_i)\big)=\mil_k^\prime\bigl(\cm_{kj}(e_j)\big)=\mil_j^\prime(e_j)$. Hence~$\factor$ is well defined, which establishes our claim.

If $\Esysprimed$ is a system in $\VS$ that is compatible with $\Edirsys$, then \cref{res:factoring_map_linear_or_vector_lattice_homomorphism} shows that the unique factoring map $\factor\colon E\to E^\prime$ as defined above in $\SET$ is, in fact, linear. We conclude that $\Esys$ is a direct limit of $\Edirsys$ in $\VS$ that is also a direct limit of the system in $\SET$.

If $\Edirsys$ is a direct system in $\VLLH$, then we proceed as above for $\VS$ to construct a direct limit $\Esys$ in $\SET$. In this case, $\widetilde E$ is a vector lattice and~$\widetilde E_0$ is an ideal in~$\widetilde E$. Hence $\widetilde E/\widetilde E_0$ is a vector lattice. The~$\maptwo_i$ are now lattice homomorphisms, and then so are the $\mil_i\colon E_i\to
\widetilde E/\widetilde E_0$. This implies that the nested union~$E$ is now a vector sublattice of $\widetilde E/\widetilde E_0$, so that the~$\mil_i$ can be viewed as lattice homomorphisms into~$E$. If $\Esysprimed$ is a direct system in $\VLLH$ that is compatible with $\Edirsys$, then  \cref{res:factoring_map_linear_or_vector_lattice_homomorphism} shows that the factoring map~$\factor$ as constructed above in $\SET$ is, in fact, a lattice homomorphism. We conclude that $\Esys$ is a direct limit of $\Edirsys$ in $\VLLH$. It is also a direct limit of the system in $\VS$.

Let $\Edirsys$ be a direct system in $\VLLH$. We claim that its direct limit $\Esys$ in $\VLLH$ as constructed above is, in addition, also a direct limit of $\Edirsys$ in the category $\POVSP$ of preordered vector spaces and positive linear maps. To see this, suppose that $\Esysprimed$ is a system in $\POVSP$ that is compatible with $\Edirsys$. Then the factoring linear map $\factor\colon E\to E^\prime$ in $\VS$ is positive. Indeed, take $e\in\pos{E}$, and choose~$i$ such that $e=\mil_i(e_i)$ for some $e_i\in E_i$. As~$\mil_i$ is a lattice homomorphism, we may suppose that $e_i\in\pos{(E_i)}$. Then $\factor(e)=\mil_i^\prime(e_i)\in\pos{(E^\prime)}$.

We now consider a direct system $\Edirsys$ in the category $\VLIPLH$ of vector lattices and interval preserving lattice homomorphisms. After the construction of the vector lattice~$E$ and the lattice homomorphisms $\mil_i\colon E_i\to E$ as for $\VLLH$, it is now natural to attempt to start with the observation that the~$\maptwo_i$ are interval preserving, and try to argue from there to show that the lattice homomorphisms~$\mil_i$ are also interval preserving. The~$\maptwo_i$, however, are interval preserving only in degenerate cases (see \cref{res:zero_connecting_morphisms_interval_preserving}). To get around this obstruction, we consider the commutative diagram
\begin{equation}\label{dia:standard_commutative_diagram}
	\begin{tikzcd}
		E_i\arrow{r}{\cm_{ji}}\arrow[two heads]{d}[swap]{\mil_i}   &  E_j\arrow[two heads]{d}{\mil_j}\\
		\mil_i(E_i)\arrow{r}[swap]{\iota_{ji}}  &   \mil_j(E_j)
	\end{tikzcd}
\end{equation}
where the bottom map is the inclusion map. \cref{res:square_lemma_vector_lattices} shows that~$\iota_{ji}$ is an interval preserving lattice homomorphism, and then \cref{res:nested_vector_lattices} implies that the inclusion maps from the $\mil_i(E_i)$ into~$E$ are also interval preserving lattice homomorphisms. We know from \cref{res:interval_preserving_maps_vector_lattices} that the surjective lattice homomorphisms $\mil_i\colon E_i\to\mil_i( E_i)$ are interval preserving, so that we can now still conclude that the compositions $\mil_i\colon E_i\to E$ are interval preserving lattice homomorphisms. If $\Esysprimed$ is a system in $\VLIPLH$ that is compatible with $\Edirsys$, then it follows from \cref{res:factoring_map_linear_or_vector_lattice_homomorphism} and \cref{res:factoring_map_interval_preserving} that the factoring map $\factor\colon E\to E^\prime$ in $\SET$ is an interval preserving lattice homomorphism.  We conclude that $\Esys$ as produced by the standard construction for direct limits of direct systems in $\VS$ is a direct limit of $\Edirsys$ in $\VLIPLH$. It is also a direct limit of the system in $\VLLH$. These two facts can already be found in \cite[Proposition~3.4]{van_amstel_van_der_walt_UNPUBLISHED:2022}.

As is easily checked, it is generally true that all direct limits of a given direct system in a category are preserved under a given functor when this holds for one particular direct limit. Using this, and the compatible isomorphisms between direct limits, the above yields the following modest improvement of some of the material in \cite{van_amstel_van_der_walt_UNPUBLISHED:2022} on direct limits.

\begin{theorem}\label{res:direct_limits_in_chain_of_vector_lattices}
	Let $\Edirsys$ be a direct system in $\bigl<$~$\VLIPLH$ / $\VLLH$ / $\VS$~$\bigr>$. Then it has a direct limit in $\bigl<$~$\VLIPLH$ / $\VLLH$ / $\VS$~$\bigr>$. If $\Esys$ is any direct limit, then $E=\bigcup_i \mil_i(E_i)$,  $\mil_i(E_i)\subseteq\mil_j(E_j)$ when $i\leq j$, and the $\mil_i(E_i)$ are $\bigl<$~ideals in~$E$ / vector sublattices of~$E$ / linear subspaces of~$E$~$\bigr>$. Furthermore, $\Esys$ is also a direct limit of $\Edirsys$ in every category to the right of $\bigl<$~$\VLIPLH$ / $\VLLH$ / $\VS$~$\bigr>$ in the chain
	\[
		\VLIPLH\hskip .3em\subset\hskip .3em\VLLH\hskip .3em\subset\hskip .3em\POVSP\hskip .3em\to\hskip .3em\VS\hskip .3em\to\hskip .3em\SET.
	\]	
\end{theorem}

It was mentioned above that the maps~$\maptwo_i$ are interval preserving only in degenerate cases. We conclude this section by tieing up this loose end.

\begin{lemma}\label{res:zero_connecting_morphisms_interval_preserving}
		Let $\Edirsys$ be a direct system in the category of vector lattices and linear maps. For every index~$i$, the following are equivalent:
		\begin{enumerate}
				\item\label{part:zero_connecting_morphisms_interval_preserving_1}
				the map $\maptwo_i\colon E_i\to\widetilde{E}$ in \cref{eq:psi_map_definition} is interval preserving;
				\item\label{part:zero_connecting_morphisms_interval_preserving_2}
				$\cm_{ji}=0$ for all $j>i$.
			\end{enumerate}
	\end{lemma}

\begin{proof}
	We prove that~\ref{part:zero_connecting_morphisms_interval_preserving_1} implies~\ref{part:zero_connecting_morphisms_interval_preserving_2}.
	Since~$\maptwo_i$ is interval preserving, so is its composition with a surjective lattice homomorphism. Hence $\xi\colon E_i\to E_i\times E_j$, defined by setting $\xi(e_i)\coloneqq(\cm_{ii}(e_i),\cm_{ji}(e_i))$, is interval preserving. Take $e_i\in\pos{E_i}$. Since $(\mil_{ii}(e_i)/2,\mil_{ji}(e_i))\in[0,\xi(e_i)]_{E_i\times E_j}$, there exists an $e_i^\prime\in[0,e_i]_{E_i}$ such that $\xi(e_i^\prime)=(\mil_{ii}(e_i)/2,\mil_{ji}(e_i))$. This implies that $\cm_{ji}(e_i/2)=\cm_{ji}(e_i)$. It follows that $\cm_{ji}=0$.
	
	It is clear that~\ref{part:zero_connecting_morphisms_interval_preserving_2} implies~\ref{part:zero_connecting_morphisms_interval_preserving_1}.
	\end{proof}

Suppose that $\Edirsys$ is a direct system in $\VLIP$ such that all maps $\maptwo_i$ are interval preserving. In view of \cref{res:zero_connecting_morphisms_interval_preserving}, there are then two possibilities:
\begin{enumerate}
	\item The index set has no largest element. Then all maps~$\mil_i^\prime$ in every compatible system $\Esysprimed$ in $\VLIP$ are the zero map.  Such a direct system has a direct limit $\Esys$ in $\VLIP$ where~$E$ is the zero space, and where all $\mil_i\colon E_i\to E$ are the zero map.
	\item The index set has a largest element~$i_{\text{l}}$. Then all maps~$\mil_i^\prime$ in every compatible system $\Esysprimed$ in $\VLIP$ are the zero map when $i\neq\mil_{i_{\text{l}}}$. Such a system has a direct limit $\Esys$ in $\VLIP$ where $E=E_{i_{\text{l}}}$, where $\mil_i\colon E_i\to E$ is the zero map when $i\neq i_{\text{l}}$, and where $\mil_{i_{\text{l}}}\colon E_{i_{\text{l}}}\to E$ is the identity map.
\end{enumerate}


\subsection{Direct limits of normed lattices}\label{subsec:direct_limits_of_normed_vector_lattices}


Analogously to the vector lattice case in \cref{subsec:direct_limits_of_vector_lattices}, direct limits can be found in a number of categories of normed lattices by using a general construction and then exploiting the additional information. The method is an adaptation of the one in \cref{subsec:direct_limits_of_vector_lattices}. In this case, the resulting standard construction produces a direct limit of a direct system $\Edirsys$ in the category $\NS$ of normed spaces and contractions that is also a direct limit in the category $\MET$ of metric spaces and contractive maps. It is as follows.

Consider the vector space
\begin{equation}\label{eq:direct_product_normed_case}
	{\widetilde E}\coloneqq\Bigl\{(e_i)\in\prod_i E_i: \sup_i\norm{e_i}<\infty\Bigr\}
\end{equation}
and supply it with the norm
\[
	\norm{(e_i)}\coloneqq \sup_{i}\norm{e_i}
\]
for $(e_i)\in {\widetilde E}$. Then
\begin{equation}\label{eq:subspace_normed_case}
	\widetilde E_0\coloneqq\left\{(e_i)\in {\widetilde E}: \lim_i\norm{e_i}=0\right\}
\end{equation}
is a closed linear subspace of~$\widetilde E$.
Since the~$\cm_{ji}$ are all contractions, we can define contractions $\maptwo_i\colon E_i\to {\widetilde E}$ by setting
\begin{equation}\label{eq:psi_map_normed_case}
	\bigl(\maptwo_i(e_i)\bigr)_j\coloneqq
	\begin{cases}
		\cm_{ji}(e_i)&\text{when }j\geq i;\\
		0&\text{else}.
	\end{cases}
\end{equation}
We let $q\colon {\widetilde E}\to {\widetilde E}/{\widetilde E}_0$ denote the quotient map between the normed spaces~$\widetilde E$ and~$\widetilde E_0$, and set $\mil_i\coloneqq q\circ\maptwo_i$. The~$\mil_i$ are contractions.
Clearly, when $e_i\in E_i$ and $e_j\in E_j$, then $\mil_i(e_i)=\mil_j(e_j)$ if and only if there exists a $k\geq i,j$ such that $\lim_{l\geq k}\norm{\cm_{li}(e_i)-\cm_{lj}(e_j)}=0$. In particular, if $\cm_{ki}(e_i)=\cm_{kj}(e_j)$ for some $k\geq i,j$, then $\mil_i(e_i)=\mil_j(e_j)$.
It is now clear that $\mil_i=\mil_j\circ \cm_{ji}$ for $i\leq j$. We define the subset
\[
	E\coloneqq\bigcup_i\mil_i(E_i)
\]
of $\widetilde E/\widetilde E_0$. As in the case of $\VS$, the fact that $\mil_i=\mil_j\circ\cm_{ji}$ for $i\leq j$ implies that the nested union~$E$ is a linear subspace of ${\widetilde E}/{\widetilde E}_0$, so that we can view the~$\mil_i$ as contractions from the~$E_i$ into the normed space~$E$.

We claim that $\Esys$ is a direct limit of $\Edirsys$ in $\MET$. The compatibility with $\Edirsys$ has again already been observed. For the universal property, suppose that the system $\Esysprimed$ in $\MET$ is compatible with $\Edirsys$. We define the (again evidently unique) factoring map $\factor\colon E\to E^\prime$ as before: take $e\in E$, choose~$i$ and $e_i\in E_i$ such that $e=\mil_i(e_i)$, and set $\factor(e)\coloneqq\mil_i^\prime(e_i)$. The requirement $\factor\circ\mil_i=\mil_i^\prime$ is met by construction, but we need to show that~$\factor$ is well defined. For this, suppose that $\mil_i(e_i)=\mil_j(e_j)$, and take $\veps>0$. There exists a $k\geq i,j$ such that $\norm{\cm_{ki}(e_i)-\cm_{kj}(e_j)}<\veps$. Since~$\mil_k^\prime$ is contractive, this implies that  $\textup{d}_{E^\prime}\bigl(\mil_i^\prime(e_i),\mil_j^\prime(e_j)\bigr)=\textup{d}_{E^\prime}\Bigl(\mil_k^\prime\bigl(\cm_{ki}(e_i)\bigr),\mil_k^\prime\bigl(\cm_{kj}(e_j)\bigr)\Bigr)<\veps$. As~$\veps$ was arbitrary, we have $\mil_i^\prime(e_i)=\mil_j^\prime(e_j)$, so that $\factor$ is well defined. It remains to be shown that~$\factor$ is contractive. Take $e_i\in E_i$ and $e_j\in E_j$, and let $\veps>0$. There exists an $\tilde e=({\tilde e}_k)\in {\widetilde E}_0$ such that
\begin{align*}
	\sup_k \lrnorm{\bigl(\maptwo_i(e_i)\bigr)_k - \bigl(\maptwo_j(e_j)\bigr)_k  +{\tilde e}_k}&=\lrnorm{\maptwo_i(e_i)-\maptwo_j(e_j) + \tilde e}\\
	&<\lrnorm{\mil_i(e_i)-\mil_j(e_j)}+\veps/2.
\end{align*}
Hence we can choose a $k\geq i,j$ such that
\[
	\lrnorm{\cm_{ki}(e_i) - \cm_{kj}(e_j)  +{\tilde e}_k}<\lrnorm{\mil_i(e_i)-\mil_j(e_j)}+\veps/2
\]
as well as $\norm{{\tilde e}_k}<\veps/2$. Then $\lrnorm{\cm_{ki}(e_i) - \cm_{kj}(e_j)}<\lrnorm{\mil_i(e_i)-\mil_j(e_j)}+\veps$, so that
\begin{align*}
	\textup{d}_{E^\prime}\Bigl(\factor\bigl(\mil_i(e_i)\bigr),\factor\bigl(\mil_j(e_j)\bigr)\Bigr)&= \textup{d}_{E^\prime}\Bigl(\factor\bigl(\mil_k(\cm_{ki}(e_i))\bigr),\factor\bigl(\mil_k(\cm_{kj}(e_j))\bigr)\Bigr)\\
	&=\textup{d}_{E^\prime}\Bigl(\mil_k^\prime(\cm_{ki}(e_i)),\mil_k^\prime(\cm_{kj}(e_j))\Bigr)\\
	&\leq \lrnorm{\cm_{ki}(e_i)-\cm_{kj}(e_j)}\\
	&<\lrnorm{\mil_i(e_i)-\mil_j(e_j)}+\veps.
\end{align*}
As~$\veps$ was arbitrary, this implies that~$\factor$ is contractive. Hence $\Esys$ is indeed a direct limit of $\Edirsys$ in $\MET$.

If $\Esysprimed$ is a system in $\NS$ that is compatible with $\Edirsys$, then \cref{res:factoring_map_linear_or_vector_lattice_homomorphism} implies that the unique factoring map $\factor\colon E\to E^\prime$ as defined above in $\MET$ is, in fact, linear. We conclude that \hskip-0.5pt$\Esys$ is a direct limit of $\Edirsys$ in $\NS$ that is also a direct limit of the system in $\MET$.\footnote{Although we shall not need it, let us still mention for the sake of completeness that it is not difficult to verify that $\norm{\mil_i(e_i)}=\limsup_{k\geq i}\norm{\cm_{ki}(e_i)}\coloneqq \inf_{j\geq i}\sup_{k\geq j}\norm{\cm_{ki}(e_i)}=\lim_{j\geq i} \sup_{k\geq j}\norm{\cm_{ki}(e_i)}$.}

If $\Edirsys$ is a direct system in $\NLLH$, then we proceed as above for $\NS$ to construct a direct limit $\Esys$ in $\MET$. In this case,~$\widetilde E$ is a normed lattice and~$\widetilde E_0$ is a closed ideal in~$\widetilde E$, so that $\widetilde E/\widetilde E_0$ is a normed vector lattice. The~$\maptwo_i$ and the~$\mil_i\colon E_i\to\widetilde E/\widetilde E_0$ are contractive lattice homomorphisms. This implies that the nested union~$E$ is a normed sublattice of $\widetilde E/\widetilde E_0$, so that we can view the~$\mil_i$ as contractive lattice homomorphisms into~$E$. If $\Esysprimed$ is a system in $\NLLH$ that is compatible with $\Edirsys$, then \cref{res:factoring_map_linear_or_vector_lattice_homomorphism} implies that the factoring map $\factor\colon E\to E^\prime$ as constructed above in $\MET$ is, in fact, a lattice homomorphism. We conclude that $\Esys$ as produced by the standard construction for direct limits of direct systems in $\NS$ is a direct limit of $\Edirsys$ in $\NLLH$. Analogously to the vector lattice case, the fact that the~$\map_i$ are lattice homomorphisms implies that it is also a direct limit in the category $\PONSP$ of preordered normed spaces and positive contractions.

Next, we consider a direct system $\Edirsys$ in the category $\NLAIPLH$ of normed lattices and contractive almost interval preserving lattice homomorphisms. Analogously to the case of vector lattices, it is only in degenerate situations that the~$\maptwo_i$ are almost interval preserving (see \cref{res:zero_connecting_morphisms_almost_interval_preserving}). In this case, the combination of diagram \eqref{dia:standard_commutative_diagram} and \cref{res:square_lemma_normed_vector_lattices} shows that the inclusion maps from the $\mil_i(E_i)$ into the $\mil_j(E_j)$ for $i\leq j$ are almost interval preserving. By \cref{res:nested_normed_vector_lattices}, the inclusion maps from the~$\mil_i(E_i)$ into~$E$ are then also almost interval preserving. Since the $\mil_i\colon E_i\to\mil_i(E_i)$ are surjective lattice homomorphisms, \cref{res:interval_preserving_maps_vector_lattices} shows that they are interval preserving. We conclude that the compositions $\mil_i\colon E_i\to E$ are almost interval preserving lattice homomorphisms. If $\Esysprimed$ is a system in $\NLAIPLH$ that is compatible with $\Edirsys$, then it follows from \cref{res:factoring_map_almost_interval_preserving} that the factoring contractive lattice homomorphism $\factor\colon E\to E^\prime$ is almost interval preserving. We thus see that the construction yields a direct limit $\Esys$ of the system in $\NLAIPLH$ that is also a direct limit in $\NLLH$.

Next, we consider a direct system $\Edirsys$ in the category $\NLIPLH$ of normed lattices and contractive interval preserving lattice homomorphisms. Then a reasoning as for $\VLIPLH$ shows that $\Esys$ as produced by the standard construction for direct limits of direct systems in $\NS$ is a direct limit of the system in $\NLIPLH$. \cref{res:factoring_map_almost_interval_preserving} implies that it is also a direct limit in $\NLAIPLH$.

The above shows that we have the following companion result to \cref{res:direct_limits_in_chain_of_vector_lattices}.

\begin{theorem}\label{res:direct_limits_in_chain_of_normed_vector_lattices}
	Let $\Edirsys$ be a direct system in $\bigl<$~$\NLIPLH$ / $\NLAIPLH$ / $\NLLH$ / $\NS$~$\bigr>$. Then it has a direct limit in  $\bigl<$~$\NLIPLH$ / $\NLAIPLH$ / $\NLLH$ / $\NS$~$\bigr>$.
	If $\Esys$ is any direct limit, then $E=\bigcup_i \mil_i(E_i)$, $\mil_i(E_i)\subseteq\mil_j(E_j)$ when $i\leq j$, and $\bigl<$~the $\mil_i(E_i)$ are ideals in~$E$ / the $\mil_i(E_i)$ are vector sublattices of~$E$ and the $\overline{\mil_i(E_i)}$ are ideals in~$E$ / the $\mil_i(E_i)$ are vector sublattices of~$E$ / the $\mil_i(E_i)$ are linear subspaces of~$E$~$\bigr>$. Furthermore, $\Esys$ is also a direct limit of $\Edirsys$ in every category to the right of $\bigl<$~$\NLIPLH$ / $\NLAIPLH$ / $\NLLH$ / $\NS$~$\bigr>$ in the chain
	\[
		\NLIPLH\hskip .3em\subset\hskip .3em\NLAIPLH\hskip .3em\subset\hskip .3em\NLLH\hskip .3em\subset\hskip .3em\PONSP\hskip .3em\to\hskip .3em\NS\hskip .3em\subset\hskip .3em\MET.
	\]	
\end{theorem}

As in \cref{subsec:direct_limits_of_vector_lattices}, we conclude by showing that the~$\maptwo_i$ are (almost) interval preserving only in degenerate cases.

\begin{lemma}\label{res:zero_connecting_morphisms_almost_interval_preserving}
	Let $\Edirsys$ be a direct system in the category of normed lattices and almost interval preserving contractions. For every index~$i$, the following are equivalent:
	\begin{enumerate}
		\item\label{part:zero_connecting_morphisms_almost_interval_preserving_1}
		the map $\maptwo_i\colon E_i\to\widetilde{E}$ in \cref{eq:psi_map_normed_case} is interval preserving;
		\item\label{part:zero_connecting_morphisms_almost_interval_preserving_2}
		the map $\maptwo_i\colon E_i\to\widetilde{E}$ in \cref{eq:psi_map_normed_case} is almost interval preserving;
		\item\label{part:zero_connecting_morphisms_almost_interval_preserving_3}
		$\cm_{ji}=0$ for all $j>i$.
	\end{enumerate}
\end{lemma}

\begin{proof}
	It is clear that~\ref{part:zero_connecting_morphisms_almost_interval_preserving_1} implies~\ref{part:zero_connecting_morphisms_almost_interval_preserving_2} and that~\ref{part:zero_connecting_morphisms_almost_interval_preserving_3} implies~\ref{part:zero_connecting_morphisms_almost_interval_preserving_1}.
	
	We prove that~\ref{part:zero_connecting_morphisms_almost_interval_preserving_2} implies~\ref{part:zero_connecting_morphisms_almost_interval_preserving_3}. Since~$\maptwo_i$ is almost interval preserving, so is a composition with a surjective continuous lattice homomorphism. Hence $\xi\colon E_i\to E_i\times E_j$, defined by setting $\xi(e_i)\coloneqq(\cm_{ii}(e_i),\cm_{ji}(e_i))$, is almost interval preserving. Take $e_i\in\pos{E_i}$.  Since $(\cm_{ii}(e_i)/2,\cm_{ji}(e_i))\in [0,\xi(e_i)]_{E_i\times E_j}$, there exists a sequence $(x_n)\subseteq[0,e_i]_{E_i}$ such that $\xi(x_n)\to (\cm_{ii}(e_i)/2,\cm_{ji}(e_i))$. That is, $x_n\to e_i/2$ and $\cm_{ji}(x_n)\to\cm_{ji}(e_i)$. This implies that $\cm_{ji}(e_i)=0$, and it follows that $\cm_{ji}=0$.	
\end{proof}

As in \cref{subsec:direct_limits_of_vector_lattices}, this leaves two possibilities for direct systems in $\NLIP$ and in $\NLAIP$ such that the~$\maptwo_i$ are (almost) interval preserving, with an analogous description of their direct limits in the pertinent category.


\subsection{Direct limits of Banach lattices}\label{subsec:direct_limits_of_Banach_lattices}


For direct limits in categories of Banach lattices, the standard construction is that of a direct limit of a system $\Edirsys$ in the category $\BS$ of Banach spaces and contractions that is also a direct limit of the system in the category $\COMMET$ of complete metric spaces and contractions.

Suppose that $\Edirsys$ is a direct system in $\BS$. In that case, one starts by carrying out the standard construction for $\NS$ that produces a direct limit in $\NS$ which is also a direct limit in $\MET$. After that, one extra step is needed because the normed space that is produced by the construction in $\NS$ need not be a Banach space. Hence it is replaced with its completion, for which there is a concrete model at hand. Its closure
\[
	E\coloneqq \overline{\bigcup_i \mil_i(E_i)}.
\]
in ${\widetilde E}/{\widetilde E}_0$ is a Banach space since~$\widetilde E$ is now a Banach space. We view the~$\mil_i$ as contractions from~$E_i$ into~$E$. Since a factoring $\bigl<$~contractive map / contraction~$\bigr>$ from $\bigcup_i\mil_i(E_i)$ into a $\bigl<$~complete metric space / Banach space~$\bigr>$ as produced by the standard construction for $\NS$ can be uniquely extended to a $\bigl<$~ contractive map / contraction~$\bigr>$ from~$E$ into that space, it is clear that $\Esys$ is a direct limit of $\Edirsys$ in $\COMMET$ as well as in $\BS$.

For a direct system $\Edirsys$ in $\BLLH$, the Banach space $\widetilde E/\widetilde E_0$ is a Banach lattice. Then so is~$E$, which is the closure of a vector sublattice. The $\mil_i\colon E_i\to E$ are contractive lattice homomorphisms. If $\Esysprimed$ is a system in $\BLLH$ that is compatible with $\Edirsys$, then \cref{res:factoring_map_linear_or_vector_lattice_homomorphism_normed_case} shows that the factoring contractive map~$\factor\colon E\to E^\prime$ in $\COMMET$ is a lattice homomorphism. We conclude that the direct limit $\Esys$ of $\Edirsys$ in $\BS$ is also a direct limit of the system in  $\BLLH$. An easy approximation argument shows that the fact that the~$\mil_i$ are lattice homomorphisms implies that it is also a direct limit in the category $\POBSP$ of preordered Banach spaces and positive contractions.

For direct systems in the category $\BLIPLH$ of Banach lattices and contractive interval preserving lattice homomorphisms, the construction for $\BS$ cannot be expected to work in general. The reason is that the inclusion map from the normed vector lattice $\bigcup_i\mil_i(E_i)$ into the Banach lattice $\overline{\bigcup_i\mil_i(E_i)}$ is only guaranteed to be almost interval preserving.

For a direct system $\Edirsys$ in the category $\BLAIPLH$ and contractive almost interval preserving lattice homomorphisms, however, the standard construction for $\BS$ \emph{does} produce a direct limit in $\BLAIPLH$, which is then also a direct limit in $\BLLH$. The argument for this is a minor modification of that for $\NLAIPLH$, using diagram \eqref{dia:standard_commutative_diagram}, \cref{res:square_lemma_normed_vector_lattices},  \cref{res:nested_normed_vector_lattices} (in a slightly different way), and \cref{res:factoring_map_almost_interval_preserving} (also in a slightly different way) again.

The above implies that we have the following companion result to  \cref{res:direct_limits_in_chain_of_vector_lattices,res:direct_limits_in_chain_of_normed_vector_lattices}.

\begin{theorem}\label{res:direct_limits_in_chain_of_banach_lattices}
	Let $\Edirsys$ be a direct system in $\bigl<$~$\BLAIPLH$ / $\BLLH$ / $\BS$~$\bigr>$. Then it has a direct limit in  $\bigl<$~$\BLAIPLH$ / $\BLLH$ / $\BS$~$\bigr>$.
	If $\Esys$ is any direct limit, then $E=\overline{\bigcup_i \mil_i(E_i)}$, $\mil_i(E_i)\subseteq\mil_j(E_j)$ when $i\leq j$, and $\bigl<$~the $\mil_i(E_i))$ are vector sublattices of~$E$ and the $\overline{\mil_i(E_i)}$ are ideals in~$E$ / the $\mil_i(E_i)$ are vector sublattices of~$E$ / the $\mil_i(E_i)$ are linear subspaces of~$E$~$\bigr>$. Furthermore, $\Esys$ is also a direct limit of $\Edirsys$ in every category to the right of $\bigl<$~$\BLAIPLH$ / $\BLLH$ / $\BS$~$\bigr>$ in the chain
	\[
		\BLAIPLH\hskip .3em\subset\hskip .3em\BLLH\hskip .3em\subset\hskip .3em\POBSP\hskip .3em\to\hskip .3em\BS\hskip .3em\subset\hskip .3em\COMMET.
	\]	
\end{theorem}


\section{Direct limits: additional results}\label{sec:direct_limits_additional_results}


\noindent The constructions in \cref{sec:direct_limits_existence_via_three_standard_constructions} do not work for the categories $\VLIP$ (the vector lattices and interval preserving linear maps), $\NLIP$ (the normed lattices and interval preserving contractions), $\NLAIP$ (the normed lattices and almost interval preserving contractions), $\BLIP$ (the Banach lattices and interval preserving contractions), $\BLAIP$ (the Banach lattices and almost interval preserving contractions), or $\BLIPLH$ (the Banach lattices and contractive interval preserving lattice homomorphisms).

For the first five of these exceptional categories, the problem is with diagram~\ref{dia:standard_commutative_diagram} to which \cref{res:square_lemma_vector_lattices} or \cref{res:square_lemma_normed_vector_lattices} are applied to see that the inclusion maps at the bottom of the diagram are (almost) interval preserving. When the connecting morphisms~$\cm_{ji}$ are not lattice homomorphisms, then the~$\maptwo_i$ need not be lattice homomorphisms, and then neither need the $\mil_i=q\circ \maptwo_i$ be. Hence the hypotheses in the corollaries need not be satisfied.

The standard construction of direct limits of direct systems of Banach lattices involves passing from $\bigcup_{i}\mil_i(E_i)$ to its closure in $\widetilde E/\widetilde E_0$. As already noted in \cref{subsec:direct_limits_of_Banach_lattices}, the corresponding inclusion map will not generally be interval preserving, but only almost interval preserving. Hence the standard construction fails (only) at the very last step and need not produce direct limits in the sixth exceptional category $\BLIPLH$.

For each of these six categories, it is an open question whether direct limits always exist. Still, for four of these categories, the results in \cref{sec:interval_preserving_and_almost_interval_preserving_linear_maps_and_lattice_homomorphisms} can be used to describe basic traits of the structure of those direct limits that \emph{do} exist.

\begin{proposition}\label{res:a_priori_structure}\quad
	\begin{enumerate}
		\item\label{part:a_priori_structure_1}
		Let $\Edirsys$ be a direct system in $\VLIP$ or $\NLIP$, and suppose that $\Esys$ is a direct limit of the system in that same category. Then $E=\bigcup_i\mil_i(E_i)$, $\mil_i(E_i)\subseteq\mil_j(E_j)$ whenever $i\leq j$, and the $\map_i(E_i)$ are ideals in~$E$.
		\item\label{part:a_priori_structure_2} Let $\Edirsys$ be a direct system in $\NLAIP$ or $\BLAIP$, and suppose that $\Esys$ is a direct limit of the system in that category. Then $E=\overline{\bigcup_i\mil_i(E_i)}$, $\mil_i(E_i)\subseteq\mil_j(E_j)$ whenever $i\leq j$, and the $\overline{\mil_i(E_i)}$ are ideals in~$E$.
	\end{enumerate}
\end{proposition}

\begin{proof}
	We prove~\ref{part:a_priori_structure_2}; the proof of~\ref{part:a_priori_structure_1} is similar but easier. We start with $\NLAIP$. Since the~$\mil_i$ are almost interval preserving, \cref{res:almost_interval_preserving_maps_normed_vector_lattices} shows that the $\overline{\mil_i(E_i)}$ are ideals in~$E$. In particular, they are vector sublattices. The compatibility of the~$\mil_i$ with the~$\cm_{ji}$ implies that $\mil_i(E_i)\subseteq \mil_j(E_j)$ whenever $i\leq j$; then also $\overline{\mil_i(E_i)}\subseteq\overline{\mil_j(E_j)}$.
	Hence the nested union $\bigcup_i\overline{\mil_i(E_i)}$ is a vector sublattice of~$E$, and then so is $\overline{\bigcup_i\overline{\mil_i(E_i)}}=\overline{\bigcup_i\mil_i(E_i)}$. We view the~$\mil_i$ as maps from the~$E_i$ into $\overline{\bigcup_i\mil_i(E_i)}$; they are then almost interval preserving by \cref{res:almost_interval_preserving_linear_maps_into_sublattices}. Suppose that $\Esysprimed$ is a system that is compatible with $\Edirsys$, and let $\factor\colon E\to E^\prime$ be the unique factoring almost interval preserving contraction. Then \cref{res:factoring_map_almost_interval_preserving} shows that the restriction of~$\factor$ to $\overline{\bigcup_i\mil_i(E_i)}$ is also an almost interval preserving contraction. By continuity and density, this restriction is the unique factoring continuous map from $\overline{\bigcup_i\mil_i(E_i)}$ into~$E^\prime$. The essential uniqueness of direct limits now implies that the inclusion map from $\overline{\bigcup_i\mil_i(E_i)}$ into~$E$ must be an isomorphism. This concludes the proof for $\NLAIP$. For $\BLAIP$ one need merely add the remark that $\overline{\bigcup_i\mil_i(E_i)}$ is a Banach space since~$E$ is.
\end{proof}

\begin{remark}\label{rem:a_priori_structure_remark}\quad
	\begin{enumerate}
		\item\label{part:a_priori_structure_remark_1}
		It does not appear to be possible to give a similar reasoning, based on the results in \cref{sec:interval_preserving_and_almost_interval_preserving_linear_maps_and_lattice_homomorphisms}, that leads to a priori knowledge about the structure of direct limits in $\BLIPLH$ or $\BLIP$.
		\item\label{part:a_priori_structure_remark_2}
		One can give similar arguments to deduce all statements about the structure of the direct limits in \cref{res:direct_limits_in_chain_of_vector_lattices}, \cref{res:direct_limits_in_chain_of_normed_vector_lattices}, and \cref{res:direct_limits_in_chain_of_banach_lattices}. The main points of these results are, therefore, the asserted existence of the direct limits and the preservation of these when passing to categories to the right in the pertinent chains.
	\end{enumerate}
\end{remark}

Some direct systems in the exceptional categories have direct limits because these coincide with direct limits in a category to the right in one of the three chains of categories. The following result is based on  \cref{res:direct_limits_in_chain_of_vector_lattices}, \cref{res:direct_limits_in_chain_of_normed_vector_lattices}, \cref{res:direct_limits_in_chain_of_banach_lattices}, \cref{res:factoring_map_interval_preserving}, and \cref{res:factoring_map_almost_interval_preserving}.

\newpage 

\begin{lemma}\label{res:limits_from_above}
	\quad
	\begin{enumerate}
		\item\label{part:limits_from_above_1}
		Let $\Edirsys$ be a direct system in $\VLIP$, and let $\Esys$ be a direct limit of this system in $\VS$. If~$E$ is a vector lattice and the $\mil_i\colon E_i\to E$ are interval preserving, then $\Esys$ is also a direct limit of $\Edirsys$ in $\VLIP$.
		\item\label{part:limits_from_above_2}
		Let $\Edirsys$ be a direct system in $\NLIP$, and let  $\Esys$ be a direct limit of this system in $\NS$. If~$E$ is a normed lattice and the $\map_i\colon E_i\to E$ are interval preserving, then $\Esys$ is also a direct limit of $\Edirsys$ in $\NLIP$.
		\item\label{part:limits_from_above_3}
		Let $\Edirsys$ be a direct system in $\NLAIP$, and let $\Esys$ be a direct limit of this system in $\NS$. If~$E$ is a normed lattice and the $\mil_i\colon E_i\to E$ are almost interval preserving, then $\Esys$ is also a direct limit of $\Edirsys$ in $\NLAIP$.
		\item\label{part:limits_from_above_4}
		Let $\Edirsys$ be a direct system in $\BLAIP$, and let  $\Esys$ be a direct limit of this system in $\BS$. If~$E$ is a Banach lattice and the $\map_i\colon E_i\to E$ are almost interval preserving, then $\Esys$ is also a direct limit of $\Edirsys$ in $\BLAIP$.
	\end{enumerate}
\end{lemma}

\begin{proof}
	We prove~\ref{part:limits_from_above_4}; the other proofs are similar. Since $\Esys$ is a direct limit of the system in $\BS$, we know from \cref{res:direct_limits_in_chain_of_banach_lattices} that $E=\overline{\bigcup_{i}\mil_i(E_i)}$. Suppose that $\Esysprimed$ is a direct system in $\BLAIP$ that is compatible with $\Edirsys$. We view $\Esysprimed$ as a direct system in $\BS$ and let $\factor\colon E\to E^\prime$ denote the unique factoring contraction in $\BS$. As the compositions $\factor\circ \mil_i=\mil_i^\prime$ are almost interval preserving, \cref{res:factoring_map_almost_interval_preserving} shows that~$\factor$ is almost interval preserving.
\end{proof}

\cref{res:direct_limits_in_chain_of_vector_lattices}, \cref{res:direct_limits_in_chain_of_normed_vector_lattices}, and \cref{res:direct_limits_in_chain_of_banach_lattices} assert that a number of direct limits of direct systems of vector lattices, normed lattices, or Banach lattices, are also direct limits of these systems in $\VS$, $\NS$, or $\BS$, respectively. In view of this, \cref{res:limits_from_above} has the following consequence.

\begin{corollary}\label{res:exceptional_direct_limits_from_below}
	\quad
\begin{enumerate}
	\item\label{part:exceptional_direct_limits_from_below_1}
	A direct limit of a direct system in $\VLIPLH$ is also a direct limit of that system in $\VLIP$.
	\item\label{part:exceptional_direct_limits_from_below_2}
	A direct limit of a direct system in $\NLIPLH$ is also a direct limit of that system in $\NLIP$.
	\item\label{part:exceptional_direct_limits_from_below_3}
	A direct limit of a direct system in $\NLIPLH$ or $\NLAIPLH$ is also a direct limit of that system in $\NLAIP$.
	\item\label{part:exceptional_direct_limits_from_below_4}
	A direct limit of a direct system in $\BLAIPLH$ is also a direct limit of that system in $\BLAIP$.
\end{enumerate}	
\end{corollary}

We conclude this section by giving four examples. They show that, for each of the five exceptional categories $\VLIP$, $\NLIP$, $\NLAIP$, $\BLIP$, and $\BLAIP$, there exists a non-trivial direct system in it for which the standard construction still produces a direct limit, even though there is no guarantee that it will. They also show that, for each of these categories, there is a direct system in it for which the standard construction does not produce a direct limit, so that the existence of these direct limits remains unclear.

\begin{example}\label{ex:construction_works_BLIP}
	Fix $1\leq p<\infty$. For the index set we take $\NN=\{1,2,\dotsc\}$, and for each $i\in\NN$ we take $E_i=\ell^p$. We let~$\mil_{ii}$ be the identity map on~$E_i$, and for $i<j$ we define $\mil_{ji}\colon E_i\to E_j$ by setting
	\begin{align}\label{eq:maps_in_ell_p}
		\begin{split}
		\map_{ji}\bigl((x_1, x_2, \dotsc)\bigr)\coloneqq\biggl( x_1, & x_2, \dotsc, x_{i-1}, \frac{x_i+x_{i+1}}{2},\frac{x_{i+2}+x_{i+3}}{2},\dotsc, \\ &\frac{x_{2j-i-2}+x_{2j-i-1}}{2}, x_{2j-i}, x_{2j-i+1}, x_{2j-i+2}, \dotsc\biggr).
	\end{split}
\end{align}
	In words:~$\mil_{ji}$ takes out the block $(x_i,\dotsc,x_{2j-i-1})$ of length $(2j-2i)$, averages the coordinates in pairs, and inserts the resulting block of the halved length $(j-i)$ again. The original $x_{2j-i}$, which is was the first coordinate after the original block, is now in position $i-1+(j-i)+1=j$. It follows from this description that $\mil_{kj}\circ\mil_{ji}=\mil_{ki}$ whenever $i\leq j\leq k$. It is clear that all~$\mil_{ji}$ are surjective interval preserving linear maps, so that we have a direct system $\Edirsys$ in $\BLIP$. However,~$\mil_{ji}$ is not a lattice homomorphism when $j>i$, so the standard construction for categories of Banach lattices is not guaranteed to work.

	The standard construction for direct limits in categories of Banach lattices starts by defining~$\widetilde E$ and~$\widetilde E_0$ as in \cref{eq:direct_product_normed_case} and \cref{eq:subspace_normed_case}, respectively. The linear map $\maptwo_i\colon E_i\to\widetilde E$ is given by \cref{eq:psi_map_normed_case}. We let $q\colon\widetilde E\to \widetilde E/\widetilde E_0$ denote the quotient map, and define $\mil_i\colon E_i\to\widetilde E/\widetilde E_0$ by setting $\mil_i\coloneqq q\circ\maptwo_i$ as in the standard construction. The candidate Banach lattice in the direct limit in $\BLIP$ is the Banach subspace $E= \overline{\bigcup_i \mil_i(E_i)}$ of $\widetilde E/\widetilde E_0$. Since $\mil_i(E_i)=\mil_j\bigl(\mil_{ji}(E_i)\bigr)=\mil_j(E_j)$ when $j\geq i$, all images $\mil_i(E_i)$ coincide. Remarkably enough, these images are a Banach sublattice of $\widetilde E/\widetilde E_0$. To see this, we introduce auxiliary maps. For $x\in\ell^p$, we define $\xi(x)\in \widetilde E$ by setting $\bigl(\xi(x)\bigr)_i\coloneqq x$ for all $i\in\NN$; this defines a lattice homomorphism $\xi\colon\ell^p\to \widetilde E$.  Next, we define $\Xi\colon\ell^p\to \widetilde E/\widetilde E_0$ by setting $\Xi\coloneqq q\circ \xi$. It is not difficult to see that~$\Xi$ is an isometric embedding of~$\ell^p$ as a Banach sublattice of $\widetilde E/\widetilde E_0$. We claim that all images $\mil_i(E_i)$ coincide with $\Xi(\ell^p)$. To see this, fix an index~$i$ and take an element $x=(x_1,x_2,\dotsc)$ of~$E_i$. Using that~$p$ is finite, it is then not too difficult to verify that $\mil_i(x)=\Xi(x^\prime)$, where $x^\prime\in\ell^p$ is given by
	\begin{equation}\label{eq:ell_p}
		x^\prime=\biggl(x_1,x_2,\dotsc, x_{i-1},\frac{x_{i}+x_{i+1}}{2},\frac{x_{i+2}+x_{i+3}}{2}, \frac{x_{i+4}+x_{i+5}}{2},\dotsc\biggr).
	\end{equation}
	Hence $\mil_i(E_i)=\Xi(\ell^p)$. Moreover, \cref{eq:ell_p} makes clear that the $\mil_i\colon E_i\to E$ are interval preserving. Suppose now that $\Edirsys$ is a direct system in $\BLIP$. Since the standard construction always produces a direct limit in $\BS$, there is a unique factoring contraction $\factor\colon E\to E^\prime$. \cref{res:factoring_map_almost_interval_preserving} then shows that~$\factor$ is interval preserving. Hence $\Esys$ as produced by the standard construction for direct systems in categories of Banach lattices is a direct limit of the direct system in $\BLIP$.

	One can also view the above direct system as a direct system in $\BLAIP$. Using that almost interval preserving contractions from the order continuous Banach lattices~$\ell^p$ into a normed lattice~$E^\prime$ are actually interval preserving (see \cref{res:duality}), we see that the standard construction for Banach lattices also produces a direct limit of the system in $\BLAIP$.

	One can also view the above direct system as a direct system in $\NLIP$. Then one does not take the closure in he standard construction, and works with $E=\bigcup_i\mil_i(E_i)$ instead. In this case, this does not change the space, and just as above for Banach lattices one see that the standard construction for normed lattices produces a direct limit in $\NLIP$, and then also in $\NLAIP$ again.
\end{example}

\begin{example}\label{ex:construction_works_VLIP}
	Consider the direct system in $\VLIP$ where the index set is~$\NN$, $E_i=c_{00}$ for all~$i$, and where the~$\cm_{ji}$ are as in \cref{ex:construction_works_BLIP}.  Arguing as in that example, the standard construction for direct systems in categories of vector lattices is seen to unexpectedly produce a direct limit in $\VLIP$ that is isomorphic to~$c_{00}$.
\end{example}

\begin{example}\label{ex:construction_fails_BLIP}
	Consider the direct system in $\BLIP$ where the index set is~$\NN$, $E_i=\ell^\infty$ for all~$i$, and where the~$\cm_{ji}$ are as in \cref{ex:construction_works_BLIP}. Again we have $E=\overline{\bigcup_i\mil_i(E_i)}=\overline{\mil_1(E_1)}$. We claim that, in this case, the Banach space~$E$ is not a vector sublattice of the Banach lattice $\widetilde E/\widetilde E_0$ in the standard construction for direct systems of Banach lattices because $\bigl|\mil_1(1,-1,1,-1,\dots)\bigr|$ is not in it. Suppose, to the contrary, that it is in~$E$. Then there exists $x\in E_1$ such that $\lrnorm{q\Bigl(\bigl|\maptwo_1(1,-1,1,-1,\dotsc)\bigr|\Bigr)-q\Bigl(\maptwo_1(x)\Bigr)}<1/2$, so that there exists $(\tilde e_j)\in \widetilde E_0$ such that $\Bigl\Vert\bigl|\maptwo_1(1,-1,1,-1,\dotsc)\bigr|-\maptwo_1(x) +\tilde e\Bigr\Vert<1/2$. Since $\lim_j\norm{\tilde e_j}=0$, there exists an~$N$ with
	\[
		\Bigl\Vert\bigl|\cm_{j1}(1,-1,1,-1,\dotsc)\bigr|-\cm_{j1}(x)\Bigr\Vert<\frac{1}{2}
	\]
	for all $j\geq N$. Now $\bigl|\cm_{j1}(1,-1,1,-1,\dotsc)\bigr|$ consists of $(j-1)$ 0's followed by 1's, so that we have that
	\[
		\lrnorm{\left(\frac{x_1+x_2}{2},\frac{x_3+x_4}{2},\dotsc,\frac{x_{2j-3}+x_{2j-2}}{2}, x_{2j-1}-1, x_{2j}-1, \dotsc\right)}<\frac{1}{2}
	\]
	for all $j\geq N$. For $j=N$ this yields that $\abs{x_{2N-1}-1}<1/2$ and $\abs{x_{2N}-1}<1/2$, so that $\abs{x_{2N-1}+x_{2N}-2}<1$. On the other hand, the fact that $\abs{x_{2j-3}+x_{2j-2}}/2<1/2$ for $j=N+1$ shows that $\abs{x_{2N-1}+x_{2N}}\!<1$, contradicting $\abs{x_{2N-1}+x_{2N}-2}\!< 1$. We conclude that~$E$ is not a vector sublattice. Hence the standard construction for categories of Banach lattices does not produce a direct limit of the system in $\BLIP$. For the same reason, it does not produce a direct limit of the system in $\BLAIP$. In addition, it does not produce a direct limit in $\NLIP$ or in $\NLAIP$: $\bigl|\mil_1(1,-1,1,-1,\dots)\bigr|$ is not in $\overline{\bigcup_i\mil_i(E_i)}$, so certainly not in $\bigcup_i\mil_i(E_i)$.	

	For each of the categories $\BLIP$, $\BLAIP$, $\NLIP$, and $\NLAIP$, it is unclear whether the above system has a direct limit in it.
\end{example}

\begin{example}\label{ex:construction_fails_VLIP}
	Consider the direct system in $\VLIP$ where the index set is~$\NN$, $E_i=\ell^\infty$ for all~$i$, and where the~$\cm_{ji}$ are as in \cref{ex:construction_works_BLIP}. Arguing as in \cref{ex:construction_fails_BLIP}, one shows that the linear subspace $E=\mil_1(E_i)$ of the vector lattice $\widetilde E/\widetilde E_0$ in the standard construction for direct systems of vector lattices is not a vector sublattice because $\bigl|\mil_1(1,-1,1,-1,\dotsc)\bigr|$ is not in it. Hence the standard construction for direct systems in categories of vector lattices does not produce a direct limit of the system in $\VLIP$. It is unclear whether the system has a direct limit in $\VLIP$.
\end{example}

	
\section{Direct limits and order continuity}\label{sec:direct_limits_and_order_continuity}


\noindent Suppose that $\Esys$ is a direct limit of a direct system $\Edirsys$ in a category of vector lattices, normed vector lattices, or Banach lattices. Which properties of the~$E_i$ are then inherited by~$E$? Not much appears to be known about this in general. It is easy to see that, for a direct system $\Edirsys$ in $\VLLH$ or $\VLIPLH$, the lattice homomorphisms $\mil_i\colon E_i\to E$ are injective if and only if all connecting lattice homomorphisms~$\cm_{ji}$ for $j\geq i$ are. In this case, $E=\bigcup_i\mil_i(E_i)$ is a union of vector sublattices that are isomorphic copies of the~$E_i$; for $\VLIPLH$ the $\mil_i(E_i)$ are even ideals in~$E$. For this situation, it is established in \cite{filter:1988} that a number of properties of the~$E_i$ are inherited by~$E$; see \cite[Theorem~3.6]{van_amstel_van_der_walt_UNPUBLISHED:2022} for an overview. The general situation where the~$\cm_{ji}$ are not necessarily injective is more demanding, however, because properties of the~$E_i$ need not be inherited by their quotients $\mil_i(E_i)$ to which results such as in \cite{filter:1988} could then be applied. We are not aware of any permanence results as in \cite{filter:1988} when the~$\mil_{ji}$ are not injective. For categories of normed lattices or Banach lattices, the situation is even a little more complicated. Here the~$\mil_i$ need not even be injective when the~$\cm_{ji}$ are, so that any results which one can derive when~$E$ is (the closure) of a union of vector sublattices or ideals do not automatically translate to the general case even with injective~$\mil_{ji}$.

In this section, we investigate the permanence of order continuity of direct limits of direct systems of order continuous Banach lattices in the general situation of not necessarily injective~$\cm_{ji}$. If $\Esys$ is a direct limit of a direct system $\Edirsys$ in a category of Banach lattices, and if the~$E_i$ are all order continuous, is~$E$ then order continuous?
	
This is not true for all categories. It can already fail in the most natural of all situations. If~$E$ is a Banach lattice such that $E=\overline{\bigcup_i E_i}$ is the closure of a union of Banach sublattices, then there is a canonical direct system in $\BLLH$ such that~$E$ is the Banach lattice in one of its direct limits. As a directed index set for the~$E_i$, we take the collection of the~$E_i$ with the ordering determined by inclusion. The connecting morphisms~$\cm_{ji}$ are then the inclusion maps from~$E_i$ into~$E_j$. If we let~$\mil_i$ denote the inclusion from~$E_i$ into~$E$, then evidently $\Esys$ is a direct limit of $\Edirsys$ in $\BLLH$. As Wickstead pointed out by the following example, it is already possible in this archetypical situation that a direct limit in $\BLLH$ of a direct system of order continuous Banach lattices is no longer order continuous.
	
\begin{example}\label{ex:wickstead}
	For $i=1,2,\dotsc$, set $c_i\coloneqq \bigl\{(x_n)\in c: x_k=x_i \text{ for all } k\geq i\bigr\}$. Then $c=\overline{\bigcup_i c_i}$. Each~$c_i$ is finite dimensional and therefore order continuous, but~$c$ is not.
\end{example}

As it turns out, the obstruction in \cref{ex:wickstead} is that the inclusion maps between the~$c_i$ are not almost interval preserving. If a direct system $\Edirsys$ in $\BLAIP$ has a direct limit $\Esys$ in $\BLAIP$ (which is not automatic), and if the~$E_i$ are all order continuous, then~$E$ is order continuous. We need a little preparation for the proof of this fact.
	
When~$E$ is a Banach lattice, we let $\gamma_E\colon E\to\nbidual{E}$ denote the canonical isometric embedding of~$E$ as a Banach sublattice of~$\nbidual{E}$. We shall use the following characterisation of order continuous Banach lattices in terms of~$\gamma_E$.
	
\begin{proposition}\label{res:characterisation_of_order_continuous_norms}
	Let~$E$ be a Banach lattice with canonical embedding $\gamma_E\colon E\to\nbidual{E}$. Then the following are equivalent:
	\begin{enumerate}
		\item\label{part:characterisation_of_order_continuous_norms_1}
		$E$ is order continuous;
		\item\label{part:characterisation_of_order_continuous_norms_2}
		$\gamma_E$ is almost interval preserving;
		\item\label{part:characterisation_of_order_continuous_norms_3}
		$\gamma_E$ is interval preserving.
	\end{enumerate}	
\end{proposition}

\begin{proof}
	The order continuity of the norm on~$E$ is equivalent to $\gamma_E(E)$ being an ideal in~$\nbidual{E}$; see \cite[Theorem~2.4]{meyer-nieberg_BANACH_LATTICES:1991}, for example. Since $\gamma_E(E)$ is closed, an appeal to the first parts of  \cref{res:interval_preserving_maps_vector_lattices,res:almost_interval_preserving_maps_and_lattice_homomorphisms_normed_vector_lattices} concludes the proof.
\end{proof}

We can now give a proof of the following result.

\begin{theorem}\label{res:order_continuity_in_BLAIP}
	Let $\Edirsys$ be a direct system in $\BLAIP$ that has a direct limit $(E,(\map_i))$ in $\BLAIP$.  If each~$E_i$ is order continuous, then so is~$E$.
\end{theorem}
	
\begin{proof}
	It follows from \cref{res:duality} that the $\nbidual{\mil_i}\colon\nbidual{E_i}\to \nbidual{E}$ are almost interval preserving contractions. The order continuity of the~$E_i$ implies that the same is true for the $\gamma_{E_i}\colon E_i\to\nbidual{E_i}$. Hence the $\nbidual{\mil_i}\circ\gamma_{E_i}\colon E_i\to \nbidual E$ are almost interval preserving contractions. The commutativity of the diagram
	\[
		\begin{tikzcd}	
			E_i\arrow{rr}{\map_{ji}}\arrow{d}[swap]{\gamma_{E_i}} &  & E_j\arrow{d}{\gamma_{E_j}}\\
			\nbidual{E_i}\arrow{rr}{\nbidual{\map_{ji}}}\arrow{rd}[swap]{\nbidual{\map_i}} &  & \nbidual{E_j}\arrow{ld}{\nbidual{\map_j}}\\
			& \nbidual{E} &
		\end{tikzcd}
	\]
	shows that the system $(\nbidual{\map_i}\circ\gamma_{E_i})$ of almost interval preserving contractions from the~$E_i$ into~$\nbidual{E}$ is compatible with $(\map_{ji})_{j\geq i}$. Hence there exists a unique almost interval preserving contraction $\factor\colon E\to\nbidual{E}$ such that $\nbidual{\map_i}\circ\gamma_{E_i}=\factor\circ\map_i$ for all~$i$. This~$\factor$ is a continuous linear map that makes the diagram
	\[
		\begin{tikzcd}
			E_i\arrow{r}{\map_i}\arrow{d}[swap]{\gamma_{E_i}} & E\arrow{d}{\factor}\\
			\nbidual{E_i}\arrow{r}{\nbidual{\map_i}} & \nbidual{E}
		\end{tikzcd}	
	\]
	commutative for all~$i$. The diagram shows that~$\factor$ agrees with~$\gamma_E$ on each $\map_i(E_i)$, and then, by continuity, also on $\overline{\bigcup_i\map_i(E_i)}$, which equals~$E$ by \cref{res:a_priori_structure}. Hence $\factor=\gamma_E$, and we see that~$\gamma_E$ is almost interval preserving. By \cref{res:characterisation_of_order_continuous_norms},~$E$ is order continuous.
\end{proof}

According to \cref{res:direct_limits_in_chain_of_banach_lattices}, a direct system in $\BLAIPLH$ has a direct limit in that category. By \cref{res:exceptional_direct_limits_from_below}, all such direct limits are also direct limits in $\BLAIP$. Hence we have the following.
	
\begin{theorem}\label{res:order_continuity_in_BLAIPLH}
	Let $\Edirsys$ be a direct system in $\BLAIPLH$. Then it has a direct limit in $\BLAIPLH$. Take a direct limit $\Esys$ in $\BLAIPLH$. If all~$E_i$ are order continuous, then so is~$E$.
\end{theorem}

We proceed to establish a few results that are particularly relevant in the context of Banach function spaces. The proof of the following preparatory result is similar to that of \cref{res:limits_from_above}.

\begin{lemma}\label{res:from_BS_down_to_AIPLH}
	Let $\Edirsys$ be a direct system in $\BLAIPLH$, and let $\Esys$ be a direct limit of this system in $\BS$. If~$E$ is a Banach lattice and the $\map_i\colon E_i\to E$ are almost interval preserving lattice homomorphisms, then $\Esys$ is also a direct limit of $\Edirsys$ in $\BLAIPLH$.
\end{lemma}
		
\begin{corollary}\label{res:order_continuity_in_archetypical_situation}
	Let~$E$ be a Banach lattice,  and let $(E_i)$ be a family of closed ideals in~$E$ such that $E=\overline{\bigcup_i E_i}$. Then the following are equivalent:
	\begin{enumerate}
		\item\label{part:order_continuity_in_archetypical_situation_1}
		all~$E_i$ are order continuous;
		\item\label{part:order_continuity_in_archetypical_situation_2}
		$E$ is order continuous
	\end{enumerate}
\end{corollary}	

\begin{proof}
	It is clear that~\ref{part:order_continuity_in_archetypical_situation_2} implies~\ref{part:order_continuity_in_archetypical_situation_1}. We prove the converse.
	According to \cref{res:almost_interval_preserving_maps_normed_vector_lattices}, the inclusion maps from the~$E_i$ into~$E$ are almost interval preserving. By \cref{res:almost_interval_preserving_linear_maps_into_sublattices}, the same is true for the inclusion maps between the~$E_i$. Hence the canonical direct system that is determined by the~$E_i$ is a system in $\BLAIPLH$. It is immediate that~$E$ and the inclusion maps from the~$E_i$ into~$E$ constitute a direct limit of this canonical system in $\BS$. An appeal to \cref{res:from_BS_down_to_AIPLH} and \cref{res:order_continuity_in_BLAIPLH} shows that~$E$ is order continuous.
\end{proof}

\cref{res:order_continuity_in_archetypical_situation} applies in particular when the~$E_i$ are projection bands.\footnote{In fact, once we know that~$E$ is order continuous, we see that the~$E_i$ \emph{are} projection bands; see \cite[Theorem~2.4.4]{meyer-nieberg_BANACH_LATTICES:1991}, for example.} We give a direct proof for this particular case, and also include a criterion for density that is of some practical relevance in Banach function spaces.	

\begin{proposition}\label{res:order_continuity_in_archetypical_situation_projection_bands}
	Let~$E$ be a Banach lattice, and let $(E_i)$ be a collection of projection bands in~$E$. The following are equivalent:
	\begin{enumerate}
		\item\label{part:order_continuity_in_archetypical_situation_projection_bands_1}
		$E=\overline{\bigcup_i E_i}$;
		\item\label{part:order_continuity_in_archetypical_situation_projection_bands_2}
		 for every $x\in E$ and every $\veps>0$, there exists an index~$i$ such that  $\norm{P_{E_i^\dc}x}<\veps$.
	\end{enumerate}
	When this is the case, the following are equivalent:
	\begin{enumerate}[label={\textup{(}\alph*\textup{)}}]
		\item\label{part:order_continuity_in_archetypical_situation_projection_bands_3}
		all~$E_i$ are order continuous;
		\item\label{part:order_continuity_in_archetypical_situation_projection_bands_4}
		$E$ is order continuous.
		\end{enumerate}
\end{proposition}

\begin{proof}
	We prove that~\ref{part:order_continuity_in_archetypical_situation_projection_bands_1} implies~\ref{part:order_continuity_in_archetypical_situation_projection_bands_2}. Take $x\in E$ and $\veps>0$. Choose an index~$i$ and a $y\in E_i$ such that $\norm{x-y}<\veps/2$. Then
	\begin{align*}
		\norm{P_{E_i^\dc}x}&=\norm{x-P_{E_i}x}\\
		&\leq \norm{x-y}+\norm{y-P_{E_i}x}\\
		&=\norm{x-y}+\norm{P_{E_i}(y-x)}\\
		&<\veps/2+\veps/2\\
		&=\veps.
	\end{align*}
	
	It is evident that~\ref{part:order_continuity_in_archetypical_situation_projection_bands_2} implies~\ref{part:order_continuity_in_archetypical_situation_projection_bands_1} because $\norm{x-P_{E_i}x}=\norm{P_{E_i^\dc}x}$ for all~$x$ and~$i$.

	We prove that~\ref{part:order_continuity_in_archetypical_situation_projection_bands_3} implies~\ref{part:order_continuity_in_archetypical_situation_projection_bands_4}. Suppose that~$E$ is not order continuous. Then, by  \cite[Theorem~2.4.2]{meyer-nieberg_BANACH_LATTICES:1991}, there exists an $x\in E$ and a disjoint sequence $(x_n)$ in~$E$ such that $0\leq x_n\leq x$ for all~$n$ and $\alpha\coloneqq \inf_n \norm{x_n}>0$. Choose an index~$i$ such that $\norm{P_{E_i^\dc}x}<\alpha/2$. Then, for all~$n$, we have $0\leq P_{E_i}x_n\leq P_{E_i}x$ and
	\(
		\norm{P_{E_i}x_n}\geq \norm{x_n}-\norm{P_{E_i^\dc}x_n}>\alpha/2.
	\)
	Again by \cite[Theorem~2.4.2]{meyer-nieberg_BANACH_LATTICES:1991}, this shows that~$E_i$ is not order continuous. This contradiction implies that~$E$ is order continuous.
	
	Since Banach sublattices of order continuous Banach lattices are order continuous (this follows from \cite[Theorem~2.4.2]{meyer-nieberg_BANACH_LATTICES:1991}, for example), it is clear that~\ref{part:order_continuity_in_archetypical_situation_projection_bands_4} implies~\ref{part:order_continuity_in_archetypical_situation_projection_bands_3}.
\end{proof}

We conclude with a result on Banach function spaces, which follows from a combination of  \cref{res:order_continuity_in_archetypical_situation} (or \cref{res:order_continuity_in_archetypical_situation_projection_bands}) and ideas from \cite{de_pagter_ricker_IN_PREPARATION} in the context of Banach function spaces over compact abelian groups. 
As usual, when $(X,\Omega,\mu)$ is a measure space, we let $\Ell^0(X,\Omega,\mu)$ denote the $\sigma$-Dedekind complete vector lattice of $\Omega$-measurable functions on~$X$, with identification of two functions when they agree~$\mu$-almost everywhere. We write~$[f]$ for the equivalence class of a measurable function~$f$. A Banach function space over $(X,\Omega,\mu)$ is an ideal in $\Ell^0(X,\Omega,\mu)$ that is supplied with a norm in which it is a Banach lattice. For a topological space~$X$, we let $\cont(X)$ resp.\ $\contc(X)$ denote the continuous resp.\ the continuous compactly supported functions on~$X$.

There are no regularity assumptions on the measure $\mu$ in the following theorem.

\begin{theorem}\label{res:application}
	Let~$X$ be a locally compact Hausdorff space such that the topologies on all its compact subsets are metrisable, let~$\Omega$ be the Borel $\sigma$-algebra of~$X$, and let $\mu\colon\Omega\to[0,\infty]$ be a measure. Let $j\colon\contc(X)\to \Ell^0(X,\Omega,\mu)$ denote the canonical map, sending~$f$ to~$[f]$. Suppose that~$E$ is a Banach function space over $(X,\Omega,\mu)$ such that $j(\contc(X))$ is contained in~$E$ and dense in~$E$. Then~$E$ is order continuous.
\end{theorem}

\begin{proof}
	For every $S\in\Omega$, we define $P_S\colon E\to E$ by setting $P_S([f])=[\chi_S f]$, where~$\chi_S$ is the indicator function of~$S$. The~$P_S$ are continuous order projections, so that their ranges $P_S(E)$ are closed ideals (even projection bands) in~$E$. When $f\in\contc(X)$, then $P_{\supp f}[f]=[f]$. Hence		
	\[
	E=\overline{j(\contc(X))}\subseteq \overline{\bigcup_{K \text{ compact }}P_K(E)}\subseteq E,
	\]
	In view of \cref{res:order_continuity_in_archetypical_situation} (or of \cref{res:order_continuity_in_archetypical_situation_projection_bands}), it will be sufficient to show that $P_K(E)$ is order continuous for every compact subset~$K$ of~$X$. Take such a~$K$. The density of $j(\contc(X))$ in~$E$ implies that $P_K(j(\contc(X)))$ is dense in $P_K(E)$. We let $\gamma_K\colon\cont(K)\to\Ell^0(X,\Omega,\mu)$ denote the map that is obtained by extending a continuous function~$f$ on~$K$ to a measurable function~$\tilde f$ on~$X$ by setting it equal to zero outside~$K$, and then defining $\gamma_K(f)\coloneqq[\tilde f]$. Since every continuous function on~$K$ can be extended to an element of $\contc(X)$ (see \cite[Theorem~20.4]{rudin_REAL_AND_COMPLEX_ANALYSIS_THIRD_EDITION:1987}, for example), we see that $P_K(j(\contc(X)))=\gamma_K(\cont(K))$. As is well known, the metrisability of~$K$ implies (and is even equivalent to) the separability of $(\cont(K),\norm{\,\cdot\,}_\infty)$; see \cite[Theorem~26.15]{jameson_TOPOLOGY_AND_NORMED_SPACES:1974}, for example. The positive map $\gamma_K\colon (\cont(K),\norm{\,\cdot\,}_\infty)\to P_K(E)$ between Banach lattices is continuous, so that $\gamma_K(\cont(K))$ is a separable subspace of $P_K(E)$. Since $\gamma_K(\cont(K))$ coincides with the dense subspace $P_K(j(\contc(X)))$ of $P_K(E)$, we conclude that $P_K(E)$ is separable. But then this $\sigma$-Dedekind complete Banach lattice $P_K(E)$ cannot contain a Banach sublattice that is isomorphic to $\ell^\infty$, so that it is order continuous by \cite[Corollary~2.4.3]{meyer-nieberg_BANACH_LATTICES:1991}.
\end{proof}


\subsection*{Acknowledgements} The authors thank Onno van Gaans and Anthony Wickstead for helpful conversations\textemdash the latter also for providing \cref{ex:wickstead}\textemdash and Ben de Pagter and Werner Ricker for making \cite{de_pagter_ricker_IN_PREPARATION} available. They thank the anonymous referee for their careful reading of the manuscript and their detailed comments.

\subsection*{Funding} During this research, the first author was supported by China Scholarship Council grant no.\  201906200101. 

\subsection*{Author Contribution} Both authors conceived the results in this paper, prepared the manuscript, and reviewed and approved the final version thereof.  

\subsection*{Conflict of Interest} No potential conflict of interests was reported by the authors.

\subsection*{Data Availability} Data sharing is not applicable to this paper as no datasets were generated or analysed during the current research.



\renewcommand{\btfs}{\mathrm}

\bibliographystyle{plain}

\urlstyle{same}

\bibliography{general_bibliography}


\end{document}